\documentclass[preprint,10pt]{elsarticle}

\usepackage[margin=1.25in]{geometry}
\usepackage{amsfonts,amssymb,amsmath}
\usepackage{amsthm}
\usepackage{bm}
\usepackage{graphicx}
\usepackage{enumerate}

\usepackage{dcolumn}
\usepackage{bm}
\usepackage{color}
\usepackage{epsfig}
\usepackage{multirow}
\usepackage{mathrsfs}
\usepackage{hyperref}
\usepackage{epstopdf}
\usepackage{subcaption}
\usepackage{autobreak}
\usepackage{tikz}
\usetikzlibrary{calc}
\usepackage{placeins}
\usepackage{float}
\usepackage{algorithm}
\usepackage{algpseudocode}
\usepackage{longtable}

\newcommand{\V}[1]{\boldsymbol{#1}}
\newcommand{\M}[1]{\mathbf{#1}}
\newcommand{\norm}[1]{\left\Vert #1\right\Vert }
\newcommand{\abs}[1]{\left|#1\right|}

\newcommand{\R}{\mathbb{R}}
\newcommand{\bigO}{{\mathcal O}}

\newtheorem{theorem}{Theorem}[section]
\newtheorem{proposition}[theorem]{Proposition}

\newtheorem{remark}[theorem]{Remark}

\newcommand{\CCM}{Center for Computational Mathematics, Flatiron Institute, New York, NY 10010, USA}
\newcommand{\CCQ}{Center for Computational Quantum Physics, Flatiron Institute, New York, NY 10010, USA}

\definecolor{DarkGreen}{RGB}{0,128,32}
\newcommand{\jk}[1]{}
\newcommand{\xg}[1]{}
\newcommand{\ab}[1]{}

\begin{document}

\begin{frontmatter}

\title{A boundary integral method for computing electron repulsion integrals in heterogeneous dielectric environments}

\author[1]{Xuanzhao Gao}
\ead{xgao@flatironinstitute.org}

\author[2]{Malte R{\"o}sner}
\ead{malte.roesner@uni-bielefeld.de}

\author[1,3]{Jason Kaye}
\ead{jkaye@flatironinstitute.org}

\author[1]{Alex Barnett}
\ead{abarnett@flatironinstitute.org}

\address[1]{\CCM}
\address[2]{Faculty of Physics, Bielefeld University, 33615 Bielefeld, Germany}
\address[3]{\CCQ}

\begin{abstract}
    We present a boundary integral method for evaluating electron repulsion integrals (ERIs) in dielectric environments of piecewise constant permittivity.
    This generalizes conventional bare-Coulomb ERIs by incorporating screening from structured dielectric surroundings, such as the substrates and/or encapsulating layers present in modern two-dimensional-material devices.
    We formulate the problem as a second-kind boundary integral equation for the induced surface charge on the dielectric interfaces, and discretize it with a Nystr\"om method on adaptively refined panels.
    The incident potential generated by each smooth orbital source density is evaluated by a Fourier truncated-kernel method in or near to its source box, and by fast multipole method (FMM) accelerated quadrature farther away.
    The linear system is solved by a Krylov iteration also accelerated by the FMM.
    We focus on cuboidal device geometries, and retain efficiency at high aspect ratios
    via close-to-touching panel quadratures.
    We further develop a distributed pipeline that computes the full ERI tensor, comprising thousands of orbital pairs, embedded in a single dielectric geometry.
    We demonstrate the method on an effective three-layer graphene slab over a Si/SiO$_2$ heterojunction, focusing on the central graphene layer represented by $198$ $p_z$-like orbitals, retaining the $2613^2$ non-negligible entries of the $198^4$-sized ERI tensor.
    The model is resolved on $1.3$ million surface unknowns, achieving a relative error of $0.004\%$, at a cost of $22$ node-hours.
\end{abstract}


\end{frontmatter}

\section{Introduction}

The four-index electron repulsion integral (ERI) tensor is a fundamental computational primitive of electronic structure theory~\cite{cramer2013essentials,helgaker2000molecular}, describing electron-electron
interactions between single-particle orbitals, i.e., the wavefunctions of individual electrons, which are spatially localized in our setting.
Given a set of $N_o$ real-valued orbitals $\V{\varphi} := \{\varphi_i\}_{i=1}^{N_o}$, the rank-four ERI tensor $\V{V} := \{V_{ijkl}\}_{i,j,k,l=1}^{N_o}$ has entries
\begin{subequations}\label{eq:four_indices}
    \begin{align}
    V_{ijkl} &= \iint_{\mathbb{R}^3 \times \mathbb{R}^3} \varphi_i(x)\varphi_j(x)\, W(x, x^\prime)\,
        \varphi_k(x^\prime)\varphi_l(x^\prime)\, \mathrm{d}x\mathrm{d}x^\prime\\
    &= \iint_{\mathbb{R}^3 \times \mathbb{R}^3} \rho_{ij}(x)\, W(x,x^\prime)\, \rho_{kl}(x^\prime)\, \mathrm{d}x\mathrm{d}x^\prime\\
    &= \int_{\mathbb{R}^3} \rho_{ij}(x)\, \phi_{kl}(x)\, \mathrm{d}x, \label{eq:target_integral}
    \end{align}
\end{subequations}
where~$W(x, x^\prime)$ is a two-body interaction kernel,
$\V{\rho} := \{\rho_{ij}\}_{i,j=1}^{N_o}$ is a set of $N_o^2$ \textit{pair densities},
each defined as
\begin{equation}\label{eq:pair_density}
    \rho_{ij}(x) = \varphi_i(x)\varphi_j(x), \quad x \in \mathbb{R}^3,
\end{equation}
and
$\V{\phi} := \{\phi_{ij}\}_{i,j=1}^{N_o}$ are the potentials generated by each pair density through the kernel $W$, thus defined as
\begin{equation}
  \phi_{ij}(x) = \int_{\mathbb{R}^3} W(x, x^\prime)\, \rho_{ij}(x^\prime)\, \mathrm{d}x^\prime.
  \label{phiij}
\end{equation}
The ERI tensor appears in many-body electronic structure methods from Hartree--Fock~\cite{roothaan1951new} and hybrid density functional theory~\cite{becke1993density} to correlated wavefunction~\cite{bartlett2007coupled}, quantum-embedding~\cite{sun2016quantum}, and lattice model~\cite{hubbard1963electron} approaches.
When the orbitals are localized, the number of non-negligible pair densities is $\mathcal{O}(N_o)$ rather than $N_o^2$, and the ERI tensor is correspondingly sparse, with $\mathcal{O}(N_o^2)$ non-negligible entries.
It can be compressed further by low-rank factorizations such as Cholesky decomposition~\cite{beebe1977simplifications}, density fitting~\cite{whitten1973coulombic,vahtras1993integral}, tensor hypercontraction~\cite{hohenstein2012tensor,lu2015compression,yeh23,yeh24,zhu2026interpolative}, or the use of structured/adaptive basis sets~\cite{white17,white19,white23,lindsey24,pan26}.

In standard settings, $W(x, x^\prime) = \frac{1}{4\pi\abs{x-x^\prime}}$ is the Coulomb kernel, corresponding to electrons interacting in vacuum or a homogeneous medium, and the evaluation of the resulting integrals is routine and highly optimized. For Gaussian basis sets, for example, they admit analytic recurrences~\cite{boys1950electronic,helgaker2000molecular}.
In this work, we consider the case in which $W$ is the Green's function for the Poisson equation in a spatially varying dielectric environment, rather than the bare Coulomb kernel.
Such a kernel is not usually available in closed form, and can be applied only by solving an electrostatic interface boundary-value problem (BVP).

Spatially varying dielectric environments arise in many experimental platforms of contemporary interest.
For example, atomically thin crystals such as graphene~\cite{novoselov2004electric} and transition metal dichalcogenide monolayers~\cite{kumar2012tunable, PhysRevLett.123.206403} are regularly exposed to oxide substrates, encapsulated in hexagonal boron nitride~\cite{dean2010boron}, and gated by nearby conducting layers.
Because a two-dimensional (2D) material screens the in-plane Coulomb interaction only weakly on its own~\cite{rytova1967screened,keldysh1979coulomb,cudazzo2011dielectric}, the effective interaction between its electrons is determined to a large degree by the screening properties of its polarizable dielectric surroundings.
This environmental screening can change quasiparticle band gaps and exciton binding energies by hundreds of meV~\cite{ugeda2014giant, raja2017coulomb}, and can be patterned laterally to define junctions within a single homogeneous monolayer, an approach known as \textit{Coulomb engineering} \cite{rosner2016two,raja2017coulomb}.
In moir\'e superlattices such as magic-angle twisted bilayer graphene \cite{cao2018unconventional}, tuning the screening environment alone can even suppress correlated insulating phases in favor of superconducting ones~\cite{stepanov2020untying, barrier2026coulomb}.

Capturing such an environment fully from first principles is infeasible, since the surrounding device is far too large to treat quantum-mechanically alongside the active material.
Practical descriptions therefore retain only a small set of localized orbitals within the ``active'' part of the system, namely the part treated quantum-mechanically, such as Wannier functions of its low-energy bands, giving an extended Hubbard-type model whose interaction coefficients are precisely the ERIs~\cite{PhysRevLett.106.236805}.
Downfolding schemes such as the constrained random phase approximation~\cite{aryasetiawan2004frequency} and quantum embedding theories for correlated defect states~\cite{sheng2022green,muechler2022quantum} account for the screening from the retained electronic degrees of freedom, while the surrounding device structure must instead be described approximately at the level of continuum electrostatics~\cite{rosner2015wannier,andersen2015dielectric}.

\begin{figure*}[htbp]
    \centering
        \includegraphics[width=0.5\textwidth]{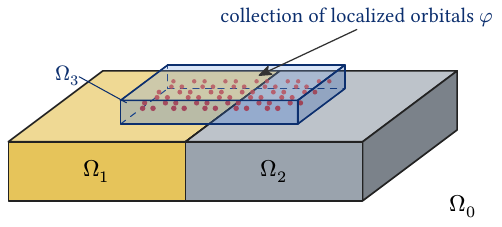}
    \caption{A typical model setup: a collection of cuboid slabs $\Omega_1$, $\Omega_2$, $\Omega_3$ with different constant dielectric permittivities, embedded in a vacuum background $\Omega_0$. The red dots mark the lattice sites of the active material, on which the localized orbitals $\V\varphi$ are centered; each orbital pair density $\rho_{ij}$ is therefore supported in a single region (here $\Omega_3$).}
    \label{fig:system}
\end{figure*}

We model the dielectric environment as a collection of (possibly-adjacent) dielectric cuboids $\V{\Omega} := \{\Omega_k\}_{k=1}^{N_k}$,
where $\Omega_k$ has constant relative permittivity $\epsilon_k$, embedded in a vacuum background
$\Omega_0 = \R^3 \setminus \bigcup_{k=1}^{N_k} \overline{\Omega_k}$ with
$\epsilon_0 = 1$. Fig.~\ref{fig:system} shows an example with $N_k=3$ materials.
This setup covers many device geometries of interest, and requires addressing high aspect ratios and edge/corner singularities,
while allowing certain other
simplifications in the boundary integral equation (BIE) implementation.
While the BIE method also applies to more general geometries with curved interfaces,
for which various accurate quadratures could be used \cite[\S9.3]{atkinson} \cite{brunoJCP01,bremer2012nystrom,GreengardONeilRachhVico2021},
the focus of the present paper is on several other challenges that arise in the ERI application.
We assume that all retained pair densities $\V{\rho}$ are supported inside a
single ``active'' region $\Omega_m$ (for example, $m=3$ in
Fig.~\ref{fig:system}). In the intended application, the basis functions
$\V{\varphi}$ are localized orbitals of this active material, obtained from an
outer electronic structure calculation such as a Wannier representation of
the low-energy bands derived from density functional theory
\cite{wannier1937structure,marzari1997maximally,marzari2012maximally}. 
Thus each pair density $\rho_{ij}$ is smooth and compactly supported inside the region. %
Dropping the $ij$ indices, a given pair density function $\rho$ induces the potential function $\phi$
which solves the Poisson BVP
\begin{equation}\label{eq:bvp}
    \begin{aligned}
        \epsilon_m \Delta \phi(x) &= -\rho(x), &&x \in \Omega_m,\\
        \Delta \phi(x) &= 0, &&x\in \mathbb{R}^3 \setminus \overline{\Omega_m},\\
        \phi\mid_{\partial \Omega_k} &= \phi\mid_{\partial \Omega_{k'}}
            &&\mbox{on }  \Gamma_{kk'},\\
        \epsilon_k \partial_{\V n} \phi\mid_{\partial \Omega_k}
            &= \epsilon_{k'} \partial_{\V n} \phi\mid_{\partial \Omega_{k'}}\quad
            &&\mbox{on }  \Gamma_{kk'},\\
        \phi(x) &\to 0, &&\abs{x}\to \infty,
    \end{aligned}
\end{equation}
where $\Gamma_{kk'} := \partial\Omega_k \cap \partial\Omega_{k'}$ is the interface between two adjacent
dielectric regions $\Omega_k$ and $\Omega_{k'}$, and $\V{n}$ denotes its normal vector
(oriented to point from $\Omega_k$ into $\Omega_{k'}$ when $k < k'$).
The first two equations are the (nondimensionalized) Poisson equation.
The third and fourth equations are continuity and flux-matching conditions.
The last line imposes zero potential at infinity, ensuring uniqueness.
Once $\phi$ has been obtained, the corresponding ERI entries can be evaluated by the inner product
\eqref{eq:target_integral} with a ``target'' pair density.

Existing continuum treatments of dielectric screening often exploit
analytically solvable models, such as effective dielectric constants,
image-charge constructions for a single planar interface, or
Rytova--Keldysh kernels for an infinite uniform slab~\cite{keldysh1979coulomb,rosner2015wannier}.
These descriptions are effective for idealized planar or otherwise separable
geometries, but the finite, multi-region, laterally-structured
configurations considered here demand a numerical solution of \eqref{eq:bvp}.

One standard approach is to discretize the Poisson equation throughout
the volume, using, for example, interface-fitted finite
differences or finite elements~\cite{leveque1994immersed,zhou2006high},
multigrid continuum-solvation solvers~\cite{andreussi2012revised}, or
particle--mesh techniques~\cite{yuan2021particle,gan2024random,gao2024broken}.
For the present problem, however, such a volumetric solve faces three difficulties:
\begin{itemize}
    \item[1.] The localized orbital densities vary on a much finer length scale
      than the surrounding dielectric structure of the device: this ratio is often $\gtrsim 10^3$.
      Any volume mesh must therefore resolve this range in each of three dimensions.
    \item[2.] The dielectric geometry contains edges and corners at which the
    electric field is singular. Resolving these local features would require
    additional 3D refinement.
    \item[3.] The problem is posed in free space with a decay condition at
      infinity; volumetric discretization thus requires artificial truncation to a sufficiently
      large domain, or specialized boundary conditions or elements for its accurate handling.
\end{itemize}

These limitations motivate a BIE approach. BIE solvers are widely used in computational electrostatics for
piecewise-homogeneous media~\cite{rokhlin1985rapid,greengard1994numerical},
with applications including Poisson--Boltzmann models of solvated
biomolecules~\cite{juffer1991electric,geng2013treecode}, polarizable continuum
models in quantum chemistry~\cite{tomasi2005quantum,cances2013domain}, and
dielectric inclusions with edges and corners~\cite{helsing2013polarizability}.
Because the permittivity is piecewise constant, the BVP \eqref{eq:bvp} can be
reformulated via Green's theorem as a BIE for an induced surface charge density defined only on
interfaces.
The reduction in dimensionality brings a large reduction in the number of unknowns---especially
those induced by multi-scale refinement---and in the cost of meshing, compared to finite elements.
The free-space condition at infinity is already built into the integral representation.
Furthermore, the BIE may be arranged to be Fredholm equation of the second kind,
so that the condition number of the resulting linear system is typically small and fixed
as the discrezation is refined.
Thus iterative solution is effective
and, although the system matrix is dense, it may be applied at each iteration
in linear cost by the fast multipole method (FMM) \cite{greengard1987fast}.

Although a BIE approach removes many of the difficulties of volumetric discretization,
some challenges remain, and motivate the techniques developed
in this work. First, the induced surface charge is singular at dielectric edges
and corners, while a localized source close to an interface produces rapidly-varying
boundary data. We handle both issues with %
adaptive
panel refinement. Second, evaluating layer potentials between nearby
panels, such as parallel facing panels in high-aspect geometries, requires {\em corrected near-singular quadratures}.
Third, each source pair density (supplied on a quadrature grid of millions of points)
generates a free-space incident potential that
must be evaluated accurately both near and far from its source box. For this, we
combine an FFT-based truncated-kernel
method (TKM)~\cite{VICO2016191} for near targets with FMM-accelerated quadrature
for far targets. Finally, computing the ERI tensor
requires many thousands of interface solves and potential evaluations on the fixed geometry.
We therefore batch the solves and reuse geometry-, source-, and target-dependent
work throughout a distributed computation pipeline.

While most of these ingredients are established, a major methodological contribution is their integration
into a scalable workflow for the automated, efficient evaluation of many dielectric-screened
ERIs in a given geometrically complex environment, to a user-prescribed accuracy.
In particular, the workflow combines adaptive interface discretization tailored to localized
pair-density sources and edge singularities,
FMM-compatible corrections for close-to-touching panels, hybrid
near/far incident-field evaluation, reuse of fixed-geometry data across many
right-hand sides, and batched distributed evaluation of the retained ERI tensor.

An outline of the paper is as follows. In Section~\ref{sec:bie}, we present the BIE formulation of the Poisson BVP \eqref{eq:bvp}, including the locally corrected Nystr\"om discretization, the adaptive panel refinement strategy, and the near-panel quadrature correction of the discrete operator. Section~\ref{sec:tkm} presents the hybrid TKM and FMM quadrature methods for evaluating the incident potential and its normal gradient.
Section~\ref{sec:assembly} explains how these tools are used in a scalable fashion:
batched interface solves over a fixed geometry, repeated potential evaluation on many target densities, and distributed computation of the sparse four-index ERI tensor.
Section~\ref{sec:numerical_results} presents numerical results:
convergence and dielectric-contrast studies on a three-material benchmark with
edges, corners, and a buried material junction, followed by computation of the
full retained ERI tensor for graphene orbital-pair densities on a
Si/SiO$_2$ heterojunction substrate.

\section{Boundary integral formulation and discretization}\label{sec:bie}

In this section, we reformulate the BVP \eqref{eq:bvp} as a second-kind BIE posed on the dielectric interfaces, and describe its discretization.
Section~\ref{sec:bie_derivation} derives the integral equation.
Section~\ref{sec:quadrature} overviews the locally corrected Nystr\"om discretization and the resulting linear system, assuming a given panel mesh.
Section~\ref{sec:adaptive_refinement} describes the adaptive construction of this panel mesh, which is guided by the orbital structure as well as the edge and corner singularities of the interface density.
Section~\ref{sec:nearfield} details the local quadrature corrections for nearby panel pairs.

\subsection{Second-kind boundary integral formulation}\label{sec:bie_derivation}

We begin by splitting the physical potential in \eqref{eq:bvp} into two parts (using language common in wave scattering theory):
\begin{equation} \label{eq:potential_split}
    \phi(x) = u_{\text{inc}}(x) + u(x).
\end{equation}
The first term is called the \textit{incident potential}, and is generated by the density $\rho$ in a homogeneous space with permittivity $\epsilon_m$,
\begin{equation}\label{eq:incident_potential_def}
    u_{\text{inc}}(x) = \frac{1}{4\pi\epsilon_m}
        \int_{\mathbb{R}^3} \frac{\rho(x')}{\abs{x-x'}}\, \mathrm{d}x',
\end{equation}
so that it satisfies
\begin{equation}\label{eq:incident_potential}
    \begin{aligned}
        \epsilon_m \Delta u_{\text{inc}}(x) &= -\rho(x), && x \in \mathbb{R}^3,\\
        u_{\text{inc}}(x) &\to 0, && \abs{x} \to \infty.
    \end{aligned}
\end{equation}
Note that the first equation holds on all of $\mathbb{R}^3$, since $\rho$ vanishes outside $\Omega_m$. Moreover, since $\rho$ is smooth and compactly supported, $u_{\text{inc}}$ is smooth throughout $\mathbb{R}^3$.
Let $\Gamma = \bigcup_{1\le k < k' \le N_k} \Gamma_{kk'}$ be the set of all interfaces.
(Note that interfaces between regions with equal permittivity are not dielectric interfaces and may be removed from $\Gamma$.)

Subtracting the conditions \eqref{eq:incident_potential} satisfied by $u_{\text{inc}}$ from the BVP \eqref{eq:bvp} for the total potential $\phi$, and using the splitting \eqref{eq:potential_split}, we see that the \textit{scattered potential} $u$ solves a Laplace BVP with inhomogeneous transmission data:
\begin{equation}\label{eq:scattered_potential}
    \begin{aligned}
        \Delta u(x) &= 0, && x \in \mathbb{R}^3 \setminus \Gamma,\\
        u^{+} - u^{-} &= 0, &&
        \mbox{on } \Gamma_{kk'},\\
        - \epsilon_{k'} \partial_{\V n} u^{+}
        + \epsilon_k \partial_{\V n} u^{-}
            &= (\epsilon_{k'} - \epsilon_k)\partial_{\V n} u_{\text{inc}},
        \quad && \mbox{on } \Gamma_{kk'},\\
        u(x) &\to 0, && \abs{x} \to \infty.
    \end{aligned}
\end{equation}
Here $u^{-}$ and $u^{+}$ denote the one-sided limits from $\Omega_k$ and $\Omega_{k'}$, respectively, consistent with the normal orientation on $\Gamma_{kk'}$.

We introduce a single layer potential representation of the scattered potential, a standard approach for solving the Laplace transmission problem with a continuous potential:
\begin{equation}\label{eq:single_layer}
    u(x) = \frac{1}{4\pi}\int_\Gamma \frac{\sigma(x')}{\abs{x-x'}}\, \mathrm{d}S(x').
\end{equation}
Here $\sigma$ is an unknown single layer density defined on $\Gamma$, and $\mathrm{d}S$ is surface area measure.
Any such representation is harmonic off $\Gamma$, continuous across $\Gamma$,
and decays at infinity, so it automatically satisfies the first, second, and fourth conditions of \eqref{eq:scattered_potential}.
For the flux-matching condition, the third line of \eqref{eq:scattered_potential},
we note that the single layer potential satisfies the standard jump relation
\begin{equation}
    \partial_{\V n} u^{\pm}(x) = \left(\mp \frac{1}{2}I + D^T\right)\sigma(x), \qquad x \in \Gamma,
\end{equation}
where $D^T$ is the on-surface adjoint double layer operator, defined via the improper integral
\begin{equation}
  (D^T \sigma)(x) = \frac{1}{4\pi} \int_\Gamma \frac{(x'-x)\cdot \V n(x)}{\abs{x-x'}^3}\sigma(x')\, \mathrm{d}S(x'),
  \qquad x\in\Gamma.
\end{equation}
For mathematical background in the Sobolev setting $\sigma\in H^{-1/2}(\Gamma)$,
sufficient for the integrable singularities arising with edges/corners of Lipschitz surfaces,
we refer the reader to, e.g., \cite[Ch.~3]{SSbook}.
Following \cite{greengard1994numerical},
substituting this jump relation into the flux-matching condition reduces it to a second-kind BIE,
\begin{equation}\label{eq:sl}
    \frac{1}{2}\gamma(x) \sigma(x) + (D^T \sigma)(x) = -\partial_{\V n} u_{\text{inc}}(x), \qquad x \in \Gamma,
\end{equation}
where~$\gamma$ is an \textit{inverse contrast function} defined on $\Gamma$ by
\begin{equation}
  \gamma(x) = \frac{\epsilon_k + \epsilon_{k'}}{\epsilon_k - \epsilon_{k'}}, \qquad x \in \Gamma_{kk'},
  \label{gammadef}
\end{equation}
obeying $|\gamma|>1$ for bounded real permittivities.
Solving \eqref{eq:sl} yields the single layer density $\sigma$,
from which $u$ is generated by \eqref{eq:single_layer},
and finally the desired total potential $\phi$ via \eqref{eq:potential_split}.

\subsection{Locally corrected Nystr\"om discretization}\label{sec:quadrature}

We discretize the integral equation \eqref{eq:sl} using a panelized Nystr\"om method.
We assume that $\Gamma$ has been partitioned into a disjoint collection $\mathcal{P}$ of $M = \abs{\mathcal{P}}$ flat rectangular panels, each living on one of the interfaces $\Gamma_{kk'}$.
The adaptive construction of this panelization will be described in Section~\ref{sec:adaptive_refinement}.
Each panel $P \in \mathcal{P}$ uses a
$p \times p$ tensor product Gauss--Legendre quadrature rule,
with nodes mapped onto the panel and weights scaled by its area.
Let
\begin{equation}
    \left\{
        x_\mu^P,\, w_\mu^P,\, \V n_\mu^P
    \right\}_{\mu=1}^{p^2}
\end{equation}
denote the mapped quadrature nodes, surface weights, and unit normals at each node on $P$.
We denote the union of all such nodes, weights, and normals as $\{x_i, w_i, \V n_i\}_{i=1}^{N}$, where $N = M p^2$.
We abbreviate
\begin{equation}
    \sigma_\mu^P = \sigma(x_\mu^P),
    \qquad
    f_\mu^P = -\partial_{\V n} u_{\text{inc}}(x_\mu^P),
\end{equation}
for the nodal values of the single-layer density and right-hand side data on panel $P$,
and collect them into vectors $\V{\sigma}_P,\, \V{f}_P \in \mathbb{R}^{p^2}$.
Stacking these vectors panel-by-panel gives
\begin{equation}
    \V{\sigma}
    =
    \begin{pmatrix}
        \V{\sigma}_{P_1} \\
        \vdots \\
        \V{\sigma}_{P_M}
    \end{pmatrix},
    \qquad
    \V{f}
    =
    \begin{pmatrix}
        \V{f}_{P_1} \\
        \vdots \\
        \V{f}_{P_M}
    \end{pmatrix}.
\end{equation}

The Nystr\"om discretization of the integral equation \eqref{eq:sl} takes the form of a block linear system
\begin{equation}\label{eq:linsys}
    \M{A}\V{\sigma} = \V{f},
\end{equation}
where $\M A \in \mathbb{R}^{N \times N}$ is composed of
square $p^2 \times p^2$ blocks $\M A_{PQ}$, each
indexed by a ``target'' panel $P$ (block row) and a ``source'' panel $Q$ (block column).
We make the decomposition
\begin{equation}\label{eq:opsplit}
    \M A = \M G + \M D^T + \M C^T,
\end{equation}
where the three terms are, respectively, the diagonal contrast term, the standard panelwise Gauss--Legendre discretization of the adjoint double layer operator, and the sparse local quadrature correction.
These terms are defined blockwise as follows:
\begin{itemize}
\item[1.] $\M G$ is block diagonal, and discretizes the first (scaled identity) term of \eqref{eq:sl}.
  The $P$th diagonal block is $\M G_{PP} = \frac{1}{2}\gamma_P \M I$, with $\gamma_P$ the contrast of the interface containing $P$.
    \item[2.] $\M D^T$ discretizes the second term of \eqref{eq:sl} using standard Gauss--Legendre quadrature:
    \begin{equation}\label{eq:block_entries}
        (\M D^T_{PQ})_{\mu\nu} = \frac{(x_\nu^Q-x_\mu^P)\cdot\V n_\mu^P}{4\pi\,|x_\mu^P-x_\nu^Q|^3}\,w_\nu^Q.
    \end{equation}
    Because all interfaces are flat in our geometry,
    $\M D^T_{PQ}\equiv0$ whenever the source and target panels are coplanar; this includes
    panels on the same interface, and in particular diagonal blocks $\M D^T_{PP}$.
    \item[3.] $\M C^T$ is a block-sparse local correction: its block $\M C^T_{PQ}$ is nonzero only when the panels $P$ and $Q$ are \emph{near} each other, or \emph{touching} at an edge or corner, in which case it corrects the standard block $\M D^T_{PQ}$ using a target-specific quadrature
    \begin{equation}\label{eq:corr_split}
        \M C^T_{PQ} = \M D^{T,\mathrm{acc}}_{PQ} - \M D^T_{PQ},
    \end{equation}
    where $\M D^{T,\mathrm{acc}}_{PQ}$ contains the accurate local quadrature weights,
    whose construction is given in Section~\ref{sec:nearfield}.
\end{itemize}
The standard rule \eqref{eq:block_entries} is high-order
accurate (with respect to panel refinement) whenever the integrand is $C^\infty$ smooth on the source panel, which holds when $P$ and $Q$ are well-separated relative to the size of $P$. This is true for all but a small fraction of panel pairs.
However, for close and touching pairs, which occur at edges, and for the large faces of thin slabs
(as in Fig.~\ref{fig:near_field}), the fixed $p \times p$ rule becomes inaccurate, and a new scheme is needed to preserve accuracy. This is the reason for introducing the correction term $\M C^T$.
This correction, and the classification of close and touching pairs, will be described in 
Section~\ref{sec:nearfield}.

With the system matrix defined, we solve the square linear system \eqref{eq:linsys} iteratively with GMRES~\cite{saad86}, using $\M{G}^{-1}$ as a diagonal
right preconditioner.
The resulting preconditioned system is
\begin{equation}\label{eq:linsys_precond}
  \M{A} \M{G}^{-1} \V{\sigma}_p = \V{f},
  \qquad \text{with}~\V{\sigma} = \M{G}^{-1} \V{\sigma}_p,
\end{equation}
so the matrix-vector multiplication at each iteration is with the matrix
\begin{equation}
    \M{A} \M{G}^{-1} = \M I + \M D^T \M G^{-1} + \M C^T \M G^{-1}.
\end{equation}
$\M G^{-1}$ is block-diagonal and can therefore be applied with $\mathcal{O}(N)$ cost.
$\M D^T$ is dense but involves source (charge)
weights that are independent of the target point,
thus can be applied to a vector at ${\cal O}(N)$ cost using a single FMM call,
with the target-normal derivative of the potential as the output.
The correction term $\M C^T$ is block sparse and can therefore be applied by a
direct sparse matrix--vector product,
with cost proportional to the number of close and touching pairs.
The latter cost is usually also linear in $N$ (since each panel has ${\cal O}(1)$ near
neighbors), and much cheaper than an FMM call.

If $\Gamma$ were smooth, $D^T$ would be compact,
so the eigenvalues of $\M A \M G^{-1}$
would cluster only at $1$,
and the GMRES convergence rate would be expected to be fixed under mesh refinement.
Corners and edges add a non-compact contribution
to $D^T$. However, its norm is known to be less than one \cite[\S6.5]{LIE}, so we still only expect
at most weak growth in iteration count.
We verify this numerically in Section~\ref{sec:convergence}.

\subsection{Adaptive panel refinement}\label{sec:adaptive_refinement}

In this section, we describe the adaptive panelization $\mathcal P$ of the set of interfaces $\Gamma$, needed to discretize the integral equation \eqref{eq:sl}.
We fix the quadrature order~$p$ and construct the panelization $\mathcal P$ using an $h$-adaptive refinement strategy, recursively subdividing each selected panel into four child panels.
Refinement is driven by two independent considerations.
First, the right-hand side~$f = -\partial_{\V n} u_{\text{inc}}$ of Eq.~\eqref{eq:sl} can vary rapidly when the source density $\rho$ is close to the interface, and the layer density $\sigma$ inherits this rapid variation.
Second, $\sigma$ develops integrable singularities near edges and corners of the dielectric interfaces.
Our corresponding refinement criteria will be data-driven and geometry-driven, respectively.

\subsubsection{Right-hand-side-driven refinement}\label{sec:rhs_refinement}

The first component of the refinement criterion is based on the local interpolation error of $f$.
On each panel $P$, let $\Pi_Pf$ be the tensor-product polynomial interpolant of degree at most $p-1$ in each local coordinate, such that
\begin{equation}\label{eq:rhs_interp}
    \Pi_Pf(x_\mu^P)=f(x_\mu^P),
    \qquad \mu=1,\ldots,p^2.
\end{equation}
We refine each panel until
\begin{equation}\label{eq:refine_criterion}
    \norm{f - \Pi_P f}_{L^2(P)} < \tau_{\mathrm{rhs}},
\end{equation}
given a prescribed local refinement threshold $\tau_{\mathrm{rhs}}$.
We estimate the norm in Eq.~\eqref{eq:refine_criterion} using an oversampled $3p \times 3p$ tensor-product Gauss--Legendre rule.

\begin{figure}[htbp]
    \centering
    \includegraphics[width=0.99\textwidth]{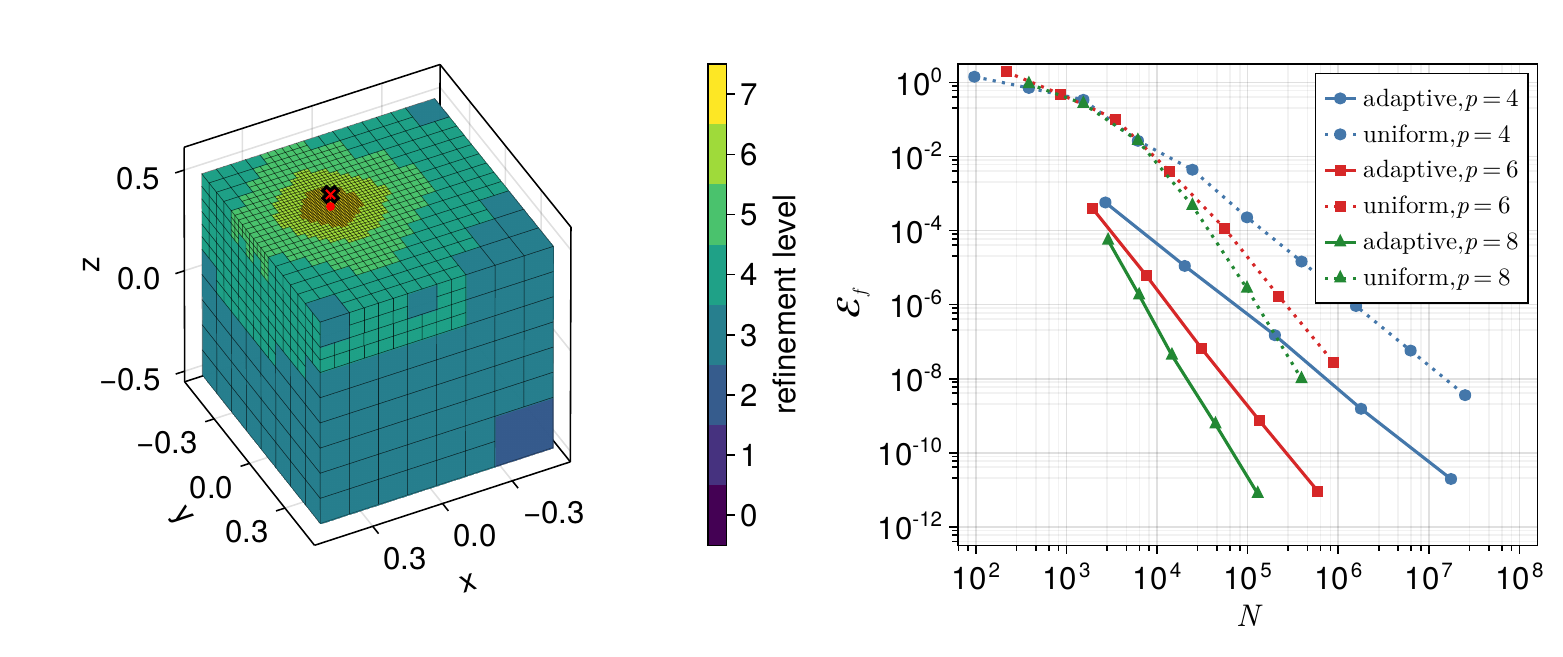}
    \caption{RHS-driven adaptive refinement for a localized source density near a dielectric interface. The dielectric is a unit cube, and the source is the normalized Gaussian $\rho(x)=(2\pi\sigma^2)^{-3/2}\exp[-\abs{x-x_0}^2/(2\sigma^2)]$, with width $\sigma = 0.03$ and center $x_0$, a height $d = 0.06$ above one face. Panels have $p\times p$ interpolation nodes. (a) An example panelization, colored by refinement level. (b) Global interpolation error \eqref{eq:Ef} versus the number of unknowns $N$ for adaptive refinement (by varying $\tau_{\mathrm{rhs}}$), compared to uniform dyadic refinement (by varying a global refinement level).
    }
    \label{fig:rhs}
\end{figure}

This adaptive refinement turns out to be crucial in our setting, as we now illustrate.
We study the absolute global interpolation $L^2$-error of the right-hand side,
\begin{equation}\label{eq:Ef}
    \mathcal{E}_f
    \;=\;
    \bigl(\sum_{P\in\mathcal{P}}\norm{f - \Pi_P f}_{L^2(P)}^2\bigr)^{1/2},
\end{equation}
for the test geometry of Fig.~\ref{fig:rhs} with data from a narrow Gaussian source $\rho$
(a simple pair density model) a typical distance above one face.
The left picture shows that the refinement concentrates near the source, as expected.
The right plot confirms that although both uniform and adaptive refinement exhibit the expected high-order convergence, %
the adaptive mesh needs $N$ about $20$ times smaller than for the uniform mesh
to reach the same error. This holds essentially independently of error level and $p$.

\subsubsection{Edge-driven refinement}\label{sec:edge_refinement}

For a rectangular panel $P$ with side lengths $s_1(P)$ and $s_2(P)$, we use
\begin{equation}
    h_P:=\max\{s_1(P),s_2(P)\}
\end{equation}
as a measure of its size.
Our geometric refinement strategy is simply to recursively subdivide
into four equal panels
any panel adjacent to an edge,
until every such panel $P$ has $h_P \le h_{\min}$,
where $h_{\min}$ is a global convergence parameter.
This results in a dyadic grading with respect to
distance from the edge (see Fig.~\ref{fig:near_field} inset),
giving ${\cal O}(h_{\min}^{-1})$ panels in total.

To predict the resulting solution error convergence rate with respect to $h_{\min}$,
we note that a dielectric edge generically induces a weak, integrable power-law singularity
in the solution density $\sigma$, matching the electric field singularity.
For simplicity, we consider the case of a single right-angle dielectric edge
(noting that the case of multi-junctions is similar \cite[(10)]{bobrovnikov73} \cite[\S5.3]{van2007electromagnetic}).
The following is derived in \ref{app:wedge} by combining classical wedge analysis with the
jump relation.

\begin{proposition}[Right-angle dielectric edge singularity]\label{prop:edge_singularity} %
  Let $d(x)$ be the distance of $x\in\Gamma_{kk'}$ from a straight right-angle edge
  (of angle either $\pi/2$ or $3\pi/2$) of the interface $\Gamma_{kk'}$
  between materials with dielectric contrast $|\gamma|>1$ defined in \eqref{gammadef}.
  Then the single-layer density obeys
\begin{equation}
    \sigma(x) = \mathcal{O}\!\left(d(x)^{\beta}\right),
    \qquad d(x)\to 0,
\end{equation}
where the singular power is
\begin{equation}
  \beta
    =\frac{2}{\pi}
    \arccos\!\left(\frac{1}{2|\gamma|}\right) - 1  \; > \; -\frac{1}{3}.
\end{equation}
\end{proposition}

We now use this for a heuristic prediction of the far-field solution error induced by
a finite $h_{\min}$.
It is known that the discretized Nystr\"om error is of the same order as the quadrature error
applied to the true solution \cite[p.~227]{LIE}.
Although panel-panel local corrections are used, the panels touching an edge
cannot accurately resolve the edge singularity.
Thus we expect the quadrature error 
to be dominated by the contribution from the unresolved innermost strip within
distance $0<r<h_{\min}$ from the edge.
By the above Proposition,
the $L^1$-norm of this contribution is bounded by
\begin{equation}\label{eq:edge_strip}
    \mathcal{E}
    \;\lesssim\;
   C \int_0^{h_{\min}} r^{\beta}\,\mathrm{d}r
    \;=\;
    \mathcal{O}\!\left(h_{\min}^{\beta+1}\right),
\end{equation}
for some unknown constant $C$, giving our convergence rate prediction for the
edge-related error.
A rigorous analysis of this would be technical
and is beyond the scope of this paper, but
the power law \eqref{eq:edge_strip} is useful in practice.

\begin{figure}[t] %
    \centering
    \includegraphics[width=0.99\textwidth]{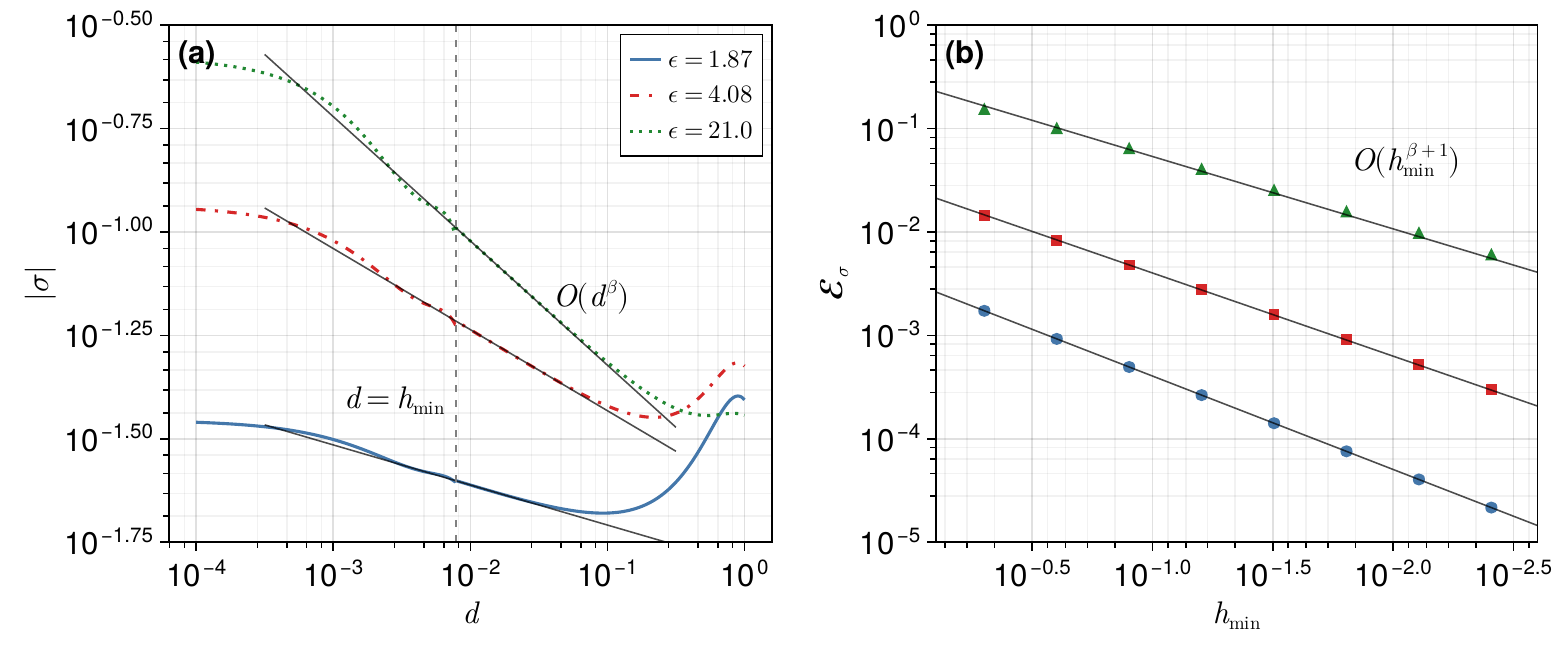}
    \caption{Edge singularity and refinement convergence test for a dielectric cube. 
    (a) Surface density magnitude $\abs{\sigma}$ along one cube face versus distance $d$ to a right-angle edge; thin black lines show the power \eqref{eq:edge_strip} and the vertical line marks $d=h_{\min}$. 
    (b) Error $\mathcal{E}_\sigma$ in the total induced surface charge versus $h_{\min}$; thin black lines show the predicted powers.}
    \label{fig:edge}
\end{figure}

To verify this, we place a unit point charge inside a dielectric cube of relative permittivity $\epsilon\in\{1.87,\,4.08,\,21.0\}$, surrounded by vacuum, so that corresponding singularity exponents are $\beta\in\{-0.1,\, -0.2,\, -0.3\}$.
We fix the quadrature order $p=4$ and vary $h_{\min}$.
By Gauss's law, the integral of surface charge equals $1-1/\epsilon$, so we choose
\begin{equation}\label{eq:Esigma}
    \mathcal{E}_\sigma
    =
    \abs{
        \int_\Gamma \sigma\,\mathrm{d}S
        -
        \left(1-\frac{1}{\epsilon}\right)
    }
\end{equation}
as an error indicator.
Figure~\ref{fig:edge}(a) shows that the density blows up with the predicted power-law $\mathcal O(d^\beta)$ over the resolved distance range, where $d$ is distance to the edge.
The singularity is stronger for larger dielectric contrast.
Figure~\ref{fig:edge}(b) shows that the predicted $\mathcal{O}(h_{\min}^{\beta+1})$
error convergence order is remarkably accurate in practice.
We have not attempted to analyze the corner singularities,
but expect their error contribution to be of lower order than the edges,
and the above results suggest this.

\subsection{Locally corrected quadrature for near-panel pairs}\label{sec:nearfield}

The adaptive mesh of Section~\ref{sec:adaptive_refinement} is designed to resolve the layer density on the $p\times p$ nodes of each panel, but the plain quadrature rule \eqref{eq:block_entries} is inaccurate for close source and target panels.
We correct the affected blocks $\M D^T_{PQ}$ as follows (note that this
does not change the number of unknowns nor their node locations).
For an ordered pair with source panel $Q$ and target panel $P$, the contribution of $Q$ at a target collocation node $x_\mu^P$ is
\[
    \int_Q K(x_\mu^P,x')\sigma(x')\,\mathrm{d}S(x'),
    \qquad \mbox{ where } \;
    K(x,x'):=\frac{(x'-x)\cdot\V n(x)}{4\pi\abs{x-x'}^3}.
\]
To express this contribution directly in terms of the unknown nodal values $\{\sigma_\nu^Q\}_{\nu=1}^{p^2}$, we interpolate
\begin{equation}\label{eq:lagrange_interp}
    \sigma_Q^{(p)}(x')
    \approx
    \sum_{\nu=1}^{p^2}\sigma_\nu^Q\,\ell_\nu^Q(x'),
    \qquad x'\in Q,
\end{equation}
where $\ell_\nu^Q$ are the product Lagrange polynomial
basis functions \cite[Ch.~5]{trefethen2019approximation} for the $p\times p$ panel nodes.
Substituting this into the panel integral defines replacement matrix elements
\begin{equation}\label{eq:accurate_block_entries}
    (\M D^{T,\mathrm{acc}}_{PQ})_{\mu\nu}
    =
    \int_Q K(x_\mu^P,x')\ell_\nu^Q(x')\,\mathrm{d}S(x'),
    \qquad \mu,\nu = 1,\dots,p^2.
\end{equation}
These integrals must now be approximated, and touching and non-touching near ordered pairs require different constructions.

For a touching pair, sources on $Q$ and targets on $P$ can approach each other closely.
We therefore compute each integral in \eqref{eq:accurate_block_entries} by adaptive 2D quadrature over $Q$ (a quad-panel version of the triangle-panel neighbor correction strategy in \cite{GreengardONeilRachhVico2021}).
The source panel $Q$ is recursively subdivided, and the same $p \times p$ tensor-product Gauss--Legendre rules are applied on its sub-panels until the estimated local
quadrature error falls below a preset parameter $\tau_{\mathrm{quad}}$.
We vectorize over the $\nu$ index for speed.

This adaptive scheme is accurate for any near pair, but expensive: each target point
in each block requires its own recursive subdivision of the source panel.
Thus, for near but non-touching pairs we apply a cheaper upsampling strategy
exploiting the fact that the kernel is now analytic on $Q$.
Let $Q_{\mathrm{up}}$ denote a $p_{\mathrm{up}}\times p_{\mathrm{up}}$ Gauss--Legendre upsampled grid on the same panel, where $p_{\mathrm{up}}>p$,
and let $\M U_{Q_{\mathrm{up}}\leftarrow Q}\in\mathbb R^{p_{\mathrm{up}}^2\times p^2}$ be
the interpolation matrix from the original $p\times p$ grid to this new grid.
Using quadrature on the upsampled grid gives
\begin{equation}\label{eq:nearfield_upsample}
    \M D^{T,\mathrm{acc}}_{PQ}
    \approx
    \M D^T_{P,Q_{\mathrm{up}}}\,
    \M U_{Q_{\mathrm{up}}\leftarrow Q},
\end{equation}
where $\M D^T_{P,Q_{\mathrm{up}}}\in\mathbb R^{p^2\times p_{\mathrm{up}}^2}$
has entries as in \eqref{eq:block_entries} but instead using the upsampled nodes and weights.

It remains to identify the non-touching ordered pairs that require correction,
then to select $p_{\mathrm{up}}$ for each.
For this, we use convergence theory involving the complex plane singularities
of the integrand (following \cite{afklinteberg2022quadrature} and other recent work),
building on analytic 1D Gauss--Legendre error estimates.
For a fixed target point $x$ and a parametrization $x'(\V\xi)$ of $Q$, with $\V\xi=(\xi_1,\xi_2)\in[-1,1]^2$, define the squared-distance function $R^2(\V\xi,x):=\abs{x'(\V\xi)-x}^2$.
For fixed $\xi_2$, its nearest complex roots in $\xi_1$ form a conjugate pair; denote the root in the upper half-plane by $t_0^{(1)}(\xi_2;x)$.
Define $t_0^{(2)}(\xi_1;x)$ analogously by fixing $\xi_1$ and solving in $\xi_2$.

For a complex singularity $t_0\notin[-1,1]$, its \emph{Bernstein parameter} is~\cite[Eq.~64]{afklinteberg2022quadrature}
\begin{equation}
    \varrho(t_0)
    =
    \left|t_0+\sqrt{t_0^2-1}\right|
    >
    1.
    \label{eq:bernstein_rho}
\end{equation}
This gives the largest Bernstein ellipse with foci at $\pm1$ within which the 1D integrand remains analytic.
The associated $p$-point Gauss--Legendre error has the exponential rate
\cite[Thm.~19.3]{trefethen2019approximation} 
\begin{equation}
    \left|I_p-I\right|
    =
    \mathcal{O}\!\left(\varrho(t_0)^{-2p}\right),
    \qquad p\to\infty.
    \label{eq:gl_rate}
\end{equation}
For the tensor-product rule, the smaller radius in the two reference directions should control the asymptotic rate.
This motivates a heuristic ``2D Bernstein'' convergence rate
(noting that more elaborate variants exist \cite{afklinteberg2022quadrature}),
\begin{equation}
    \varrho(x,Q)
    :=
    \min\left\{
        \min_{\xi_2\in[-1,1]}
        \varrho\!\left(t_0^{(1)}(\xi_2;x)\right),
        \;
        \min_{\xi_1\in[-1,1]}
        \varrho\!\left(t_0^{(2)}(\xi_1;x)\right)
    \right\}.
    \label{eq:rho_target_panel}
\end{equation}
As the target approaches $Q$, then $\varrho(x,Q)\to 1^+$, making the convergence arbitrarily
slow.
Conversely, a larger $\rho$ corresponds to a larger analyticity neighborhood and faster Gauss--Legendre convergence.

\begin{figure}[t] %
    \centering
    \includegraphics[width=0.99\textwidth]{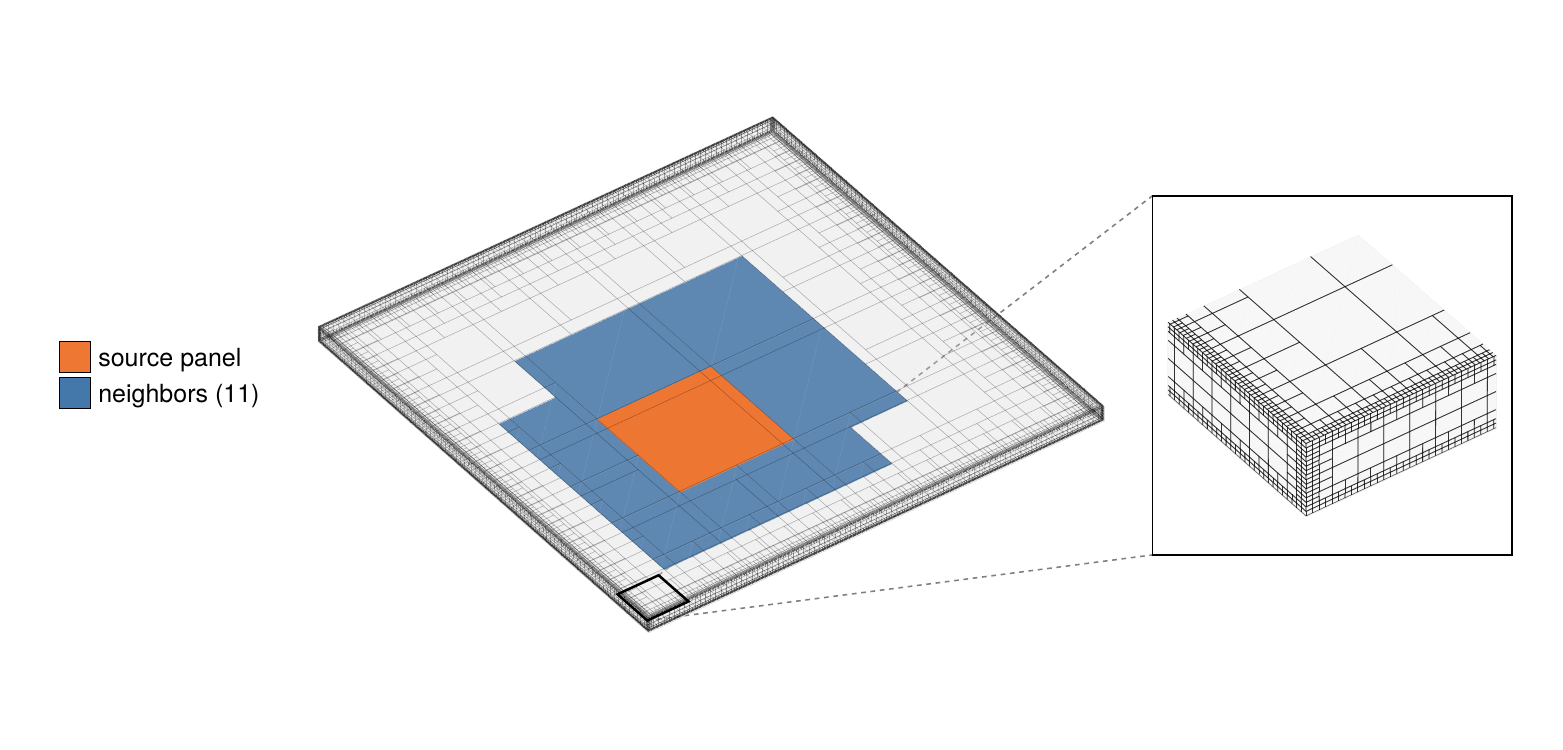}
    \vspace{-8ex}
    \caption{The nearness criterion for one fixed source panel $Q$ (orange). Target panels $P$ satisfying \eqref{eq:bernstein_criterion} are shown in blue. The cuboid has dimensions $30\times30\times1$, the Gauss--Legendre order is $p=8$, and $\tau_{\mathrm{quad}}=10^{-6}$, giving $\varrho_\star=\tau_{\mathrm{quad}}^{-1/(2p)}\approx2.37$. 
      The zoomed inset shows the dyadic edge refinement
      down to panel size $h_{\min}$.
    }
    \label{fig:near_field}
\end{figure}

Taking the worse-case Bernstein parameter over the target points in panel $P$ leads to
\begin{equation}
    \varrho_{\min}(P,Q)
    :=
    \min_{\mu=1,\ldots,p^2}\varrho(x_\mu^P,Q),
    \label{eq:rho_min}
\end{equation}
with expected quadrature error $\mathcal{O}\!\left(\varrho_{\min}(P,Q)^{-2p}\right)$ for a
$p\times p$ rule.
Equating this error with $\tau_{\mathrm{quad}}$, ignoring the prefactor,
motivates our classifying a non-touching ordered pair $(P,Q)$ as \emph{near} when
\begin{equation}\label{eq:bernstein_criterion}
    \varrho_{\min}(P,Q)\leq\varrho_\star,
    \qquad
    \varrho_\star=\tau_{\mathrm{quad}}^{-1/(2p)}.
\end{equation}
Note that since $\varrho_{\min}(P,Q)$ need not equal $\varrho_{\min}(Q,P)$,
\emph{nearness} is not a symmetric relation.
For each near ordered pair $(P,Q)$ we then choose the upsampled order by inverting the asymptotic rate formula:
\begin{equation}\label{eq:pup}
    p_{\mathrm{up}}
    =
    \left\lceil
    \frac{-\log\tau_{\mathrm{quad}}}
         {2\log\varrho_{\min}(P,Q)}
    \right\rceil.
\end{equation}
With $p_{\mathrm{up}}$ chosen, $\M D_{PQ}^{T,\mathrm{acc}}$ is filled via
\eqref{eq:nearfield_upsample}, then $\M C^T_{PQ}$ is filled via \eqref{eq:corr_split}.
For an efficient implementation, a k-d tree on panel centroids is used to generate spatially-nearby candidate pairs.
For each ordered panel pair $(P,Q)$ and each target node $x_\mu^P$, we approximate the minimization in Eq.~\eqref{eq:rho_target_panel} over the $p$ Gauss--Legendre nodes in each direction on $Q$. 
We then apply Eq.~\eqref{eq:rho_min} and the criterion in Eq.~\eqref{eq:bernstein_criterion}.
Figure~\ref{fig:near_field} highlights the resulting set of target panels $P$ that require
near corrections for one large source panel $Q$ in a thin cuboid geometry.

\begin{figure}[t]  %
    \centering
    \includegraphics[width=0.95\textwidth]{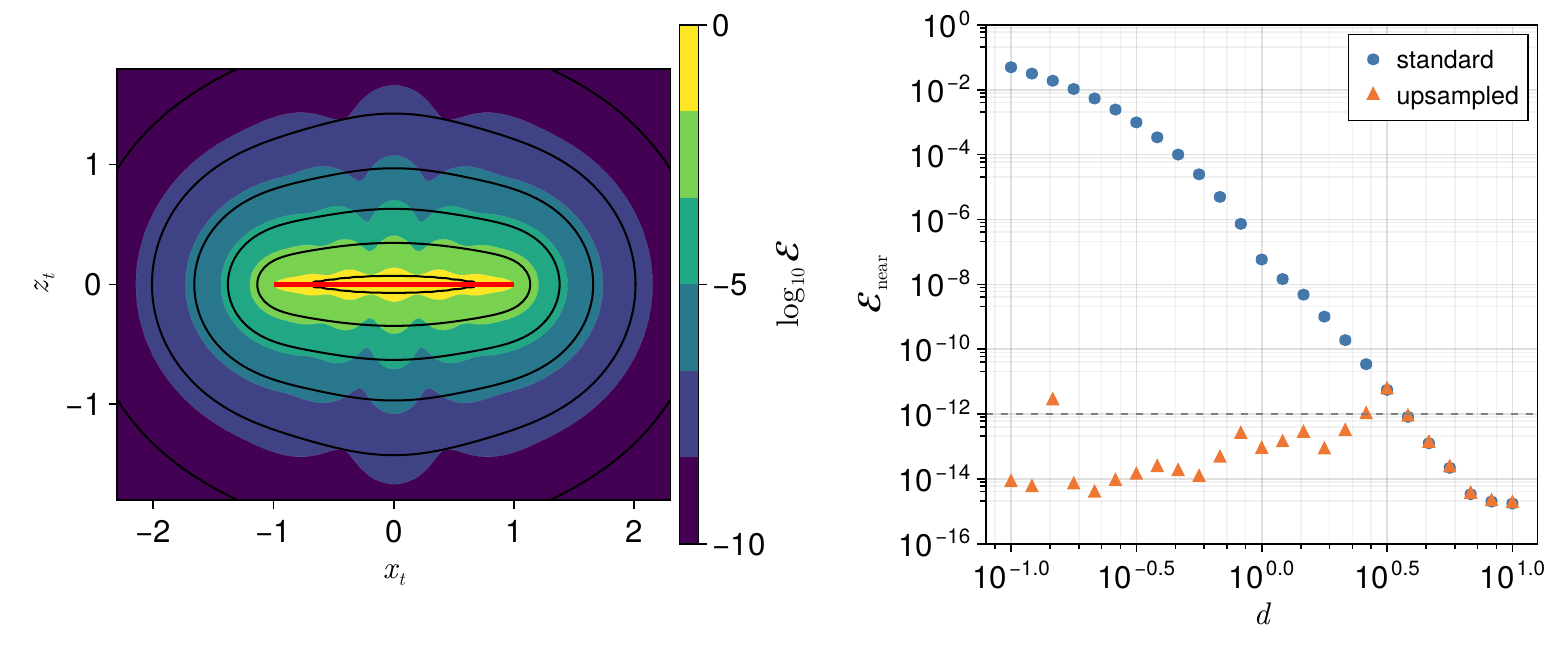}
    \caption{Test of the 2D Bernstein error prediction and non-touching quadrature correction in Section~\ref{sec:nearfield}.
      The source panel is $Q=[-1,1]^2\times\{0\}$ with density $\sigma_Q(\xi_1,\xi_2)=e^{-\xi_1^2-\xi_2^2}$ and $p=8$. 
      (a) Colors show pointwise error $\log_{10}(|I_p-I_{\mathrm{ref}}|/\max|I_{\mathrm{ref}}|)$ of the standard Gauss--Legendre rule for the single-layer kernel $1/(4\pi\abs{x-x'})$ on the plane $x_2=0$ slicing through $Q$;
      black contours show, at the same levels,
      the prediction $\varrho(x,Q)^{-2p}$ using \eqref{eq:rho_target_panel}.
      (b) Relative $L^2$ errors \eqref{eq:E_near} 
      in the adjoint double-layer block, for the standard rule (circles) and the upsampled rule \eqref{eq:pup} (triangles), as a function of panel separation $d$ (see text).
      Errors were computed using an adaptive cubature reference with 14-digit relative tolerance~\cite{HCubature}.
    }
    \label{fig:near_panel_pair}
\end{figure}

We now verify the above error estimates and the resulting near quadrature correction.
Figure~\ref{fig:near_panel_pair}(a)
provides evidence that the above 2D Bernstein radius gives a reasonable prediction of the
error using the standard (not upsampled) $p\times p$ rule.
Their agreement over the plotted range shows that the Bernstein radius captures the observed spatial error landscape.
Figure~\ref{fig:near_panel_pair}(b) validates the non-touching-panel upsampling strategy,
for parallel panels separated by a range of gaps $d$:
\begin{equation}
    Q = [-1,1]^2 \times \{0\},
    \qquad
    P = [-1,1]^2 \times \{d\}.
\end{equation}
We plot the relative $L^2$ error
\begin{equation}\label{eq:E_near}
    \mathcal{E}_{\mathrm{near}}(d)
    \;=\;
    \frac{\norm{\bigl(\M D^{T,\bullet}_{PQ} - \M D^{T,\mathrm{ref}}_{PQ}\bigr)\V\sigma_Q}_2}
         {\norm{\M D^{T,\mathrm{ref}}_{PQ}\V\sigma_Q}_2}.
\end{equation}
Here $\M D^{T,\bullet}_{PQ}$ is either the standard $p=8$ block or the upsampled replacement,
choosing $\tau_{\mathrm{quad}}=10^{-12}$.
The uncorrected rule loses accuracy as $d$ shrinks, whereas the upsampled rule maintains accuracy throughout the tested non-touching regime.
Over the same range, $p_{\mathrm{up}}$ increases from $16$ at $d=1$ to $139$ at $d=0.1$.

\section{Evaluation of incident potential and its normal gradient}
\label{sec:tkm}

This section describes the evaluation of the incident
potential Eq.~\eqref{eq:incident_potential_def} and its gradient, required respectively
for the evaluation of the physical
potential \eqref{eq:potential_split} and the right-hand side of the BIE \eqref{eq:sl}.
This is done most efficiently by separation into targets \emph{near} to the pair density support, handled by the truncated-kernel method (TKM), and the remaining \emph{far} targets, handled by
particle-FMM accelerated quadrature.

\subsection{Far-field quadrature evaluation and near/far splitting}

We first describe a trapezoidal quadrature approximation
for the incident potential Eq.~\eqref{eq:incident_potential_def} and its gradient at far-field
targets.
Let $\rho$ be a smooth source density whose support is contained in a finite cuboidal box $\mathcal{B}$.
Choosing the center of $\mathcal B$ as the origin, we write
\begin{equation} \label{eq:source_box}
    \mathcal B
    =[-l_x/2,l_x/2] \times [-l_y/2,l_y/2] \times [-l_z/2,l_z/2].
\end{equation}
We sample $\rho$ on a uniform Cartesian grid in $\mathcal B$ with spacing $h$ in all three directions.
The numbers of grid points in each direction are $N_x$, $N_y$, and $N_z$, respectively.
The grid nodes are
\begin{equation}\label{eq:source_grid}
  \textstyle
    y_{a, b, c}
    =\bigl(-\frac{l_x}{2}+ a h,\,-\frac{l_y}{2}+b h,\,-\frac{l_z}{2}+c h\bigr),
    \quad
    a = 0,1,\dots,N_x-1,\;
    b = 0,1,\dots,N_y-1,\;
    c = 0,1,\dots,N_z-1.   
\end{equation}
We enumerate these nodes as $\{y_j\}_{j=1}^{N_\rho}$, where $N_\rho=N_xN_yN_z$.
Let $\rho_h$ denote the trigonometric interpolant~\cite[\S 9]{doi:10.1137/1.9781611971514} constructed from the grid samples of $\rho$ with spacing $h$ (since $\rho$ decays smoothly to zero, it may be treated as periodic). Then we choose the grid spacing $h$ small enough so that
\begin{equation}
    \frac{\norm{\rho_{h/2} - \rho_h}_{L^2(\mathcal B)}}{\norm{\rho_{h/2}}_{L^2(\mathcal B)}} < \tau_\rho,
\end{equation}
where $\tau_\rho$ is a prescribed tolerance. This is equivalent to insuring that the tails of the Fourier
transform of $\rho$ are small beyond the Nyquist frequency magnitude $\pi/h$.

For far-field targets, we simply apply the tensor-product trapezoidal rule, discretizing the volume integral \eqref{eq:incident_potential_def} on $\{y_j\}_{j=1}^{N_\rho}$.
All weights are $\omega_\rho=h^3$, giving
\begin{equation}
    u_{\mathrm{inc}}(x)
    \approx
    \frac{\omega_\rho}{4 \pi \epsilon_m}
    \sum_{j=1}^{N_\rho} \rho(y_j) \frac{1}{\abs{x-y_j}},
    \qquad 
    \nabla u_{\mathrm{inc}}(x)
    \approx
    -\frac{\omega_\rho}{4\pi\epsilon_m}
    \sum_{j=1}^{N_\rho}
    \rho(y_j)\frac{x-y_j}{\abs{x-y_j}^3},
    \label{eq:particle_gradient}
\end{equation}
which together are evaluated efficiently at all targets $x$ with a single FMM call.

For targets far from the support of $\rho$, the integrands are smooth, so this trapezoidal rule converges superalgebraically~\cite[\S 10]{doi:10.1137/130932132}.
However, this convergence is not uniform as the target approaches the support of $\rho$, because the integrands are nearly singular and become underresolved on the fixed $h$ grid.
To distinguish between far and near targets, we define a near-field region $\mathcal B_{\mathrm{pad}}$ by padding $\mathcal B$ by $h_{\mathrm{pad}} = c_{\mathrm{pad}} h$ in each of the six directions.
Motivated by the analogous trapezoidal-rule error analysis for layer
potentials~\cite{barnett2014evaluation}, which proves that the error decays to machine precision
as the target moves several grid spacings from the source, we take
$c_{\mathrm{pad}}=5$ as a practical choice.
For the set of targets inside $\mathcal B_{\mathrm{pad}}$, we instead apply the TKM to evaluate the volume integral.

\subsection{Truncated kernel method for near-field volume integral evaluation}\label{sec:tkm_near_field}

We introduce a Fourier-based truncated kernel method, inspired by Vico et al.~\cite{VICO2016191}, to evaluate the volume integral in \eqref{eq:incident_potential_def} for targets inside $\mathcal B_{\mathrm{pad}}$.
By the convolution theorem, this volume integral has the Fourier representation
\begin{equation}\label{eq:incident_potential_fourier}
    u_{\mathrm{inc}}(x) = \frac{1}{(2\pi)^3\epsilon_m}
    \int_{\mathbb R^3}
    \frac{\widehat\rho(k)}{\abs{k}^2}e^{ik\cdot x}\,\mathrm{d}k,
\end{equation}
where $\abs{k}^{-2}$ is the Fourier transform of the Coulomb kernel, and
the source Fourier transform $\widehat\rho$ is
\begin{equation}
    \widehat\rho(k)
    =
    \int_{\mathcal B}
    \rho(x')e^{-ik\cdot x'}\,\mathrm{d}x'.
\end{equation}
This representation avoids direct integration against the singular real-space kernel at $x'=x$.
Since $\rho$ is smooth and compactly supported, $\widehat\rho$ decays rapidly so that the Fourier integral can be truncated to a bounded frequency domain, which is a cube of side length $2 k_{\max}$ centered at the origin.
A uniform-grid trapezoidal discretization then permits FFT-based evaluation.
Two issues need to be handled: \eqref{eq:incident_potential_fourier} has a singularity at $k=\bf{0}$,
and the aliasing error induced by Fourier-space sampling must be controlled.

To remove the Fourier-space singularity without changing the potential in $\mathcal B_{\mathrm{pad}}$, 
we truncate the Coulomb kernel in real space beyond the largest distance from a source $x'\in\mathcal B$ to a target $x\in\mathcal B_{\mathrm{pad}}$,
\begin{equation}
\label{eq:L_padded}
    L
    =
    \sqrt{
        (l_x+h_{\mathrm{pad}})^2+
        (l_y+h_{\mathrm{pad}})^2+
        (l_z+h_{\mathrm{pad}})^2
    }.
\end{equation}
The radially truncated Coulomb kernel is then defined by
\begin{equation}
    G_L(x)
    =
    \begin{cases}
        \frac{1}{4 \pi \abs{x}}, & \abs{x}\leq L,\\
        0, & \abs{x}>L.
    \end{cases}
\end{equation}
We then define
\begin{equation}\label{eq:tkm_convolution}
    u_L(x) := \frac{1}{\epsilon_m} \int_{\mathcal B} G_L(x-x')\rho(x')\,\mathrm{d}x' 
            = \frac{1}{(2\pi)^3\epsilon_m} \int_{\mathbb R^3} \widehat G_L(k) \widehat\rho(k)e^{ik\cdot x}\,\mathrm{d}k.
\end{equation}
Since for every $x'\in\mathcal B$ and $x\in\mathcal B_{\mathrm{pad}}$, $G_L(x-x')=G(x-x')$, 
$u_{\mathrm{inc}}(x) = u_L(x)$ for all $x\in\mathcal B_{\mathrm{pad}}$.
The Fourier transform of $G_L$ is given by \cite{VICO2016191}
\begin{equation}
\label{eq:G_L_hat}
    \widehat G_L(k)
    =
    \begin{cases}
        L^2/2, & k=0,\\
        2\left(\dfrac{\sin(L\abs{k}/2)}{\abs{k}}\right)^2,
        & k\neq0,
    \end{cases}
\end{equation}
which is now smooth.
We discretize the Fourier integral \eqref{eq:tkm_convolution} using a truncated trapezoidal rule on a Cartesian grid with spacing $\Delta k_\alpha$ in each direction (again here $a,b,c\in\mathbb{Z}$),
\begin{equation}
\label{eq:fourier_lattice}
    k_{a, b, c}
    =
    (a \Delta k_x, b \Delta k_y, c \Delta k_z),
    \qquad
    -M_x \le a \le M_x,~ -M_y \le b \le M_y,~ -M_z \le c \le M_z,
\end{equation}
where~$M_\alpha = \lceil \frac{k_{\max}}{\Delta k_\alpha} \rceil$ for $\alpha \in \{x, y, z\}$.
We explain the optimal selection of the Fourier grid spacings $\Delta k_\alpha$ in the next section.
We enumerate this Fourier grid as $\{k_\mu \}_{\mu = 1}^{M_k}$, where $M_k = (2 M_x + 1)(2 M_y + 1)(2 M_z + 1)$.
At each Fourier grid point, we approximate the Fourier transform of the source density using its existing grid \eqref{eq:source_grid}:
\begin{equation}
\label{eq:rho_hat_discrete}
    \widehat\rho(k_\mu)
    \approx
    \widehat \rho_h(k_\mu)
    =
    \omega_\rho
    \sum_{j=1}^{N_\rho}
    \rho(y_j)e^{-ik_\mu\cdot y_j},
    \qquad \mu = 1,\dots,M_k,
\end{equation}
which can be evaluated via a single padded 3D FFT.
Using this discretized $\hat\rho$ grid,
the approximation to \eqref{eq:tkm_convolution}
gives the incident potential
\begin{equation}\label{eq:tkm_u_inc_quadrature}
    u_{\mathrm{inc}}(x)
    \approx
    \frac{\Delta k_x\Delta k_y\Delta k_z}{(2\pi)^3\epsilon_m}
    \sum_{\mu = 1}^{M_k}
    \widehat G_L(k_\mu)\widehat\rho_h(k_\mu)e^{ik_\mu\cdot x}.
\end{equation}
Analytic differentiation of \eqref{eq:tkm_u_inc_quadrature} gives the corresponding gradient approximation:
\begin{equation}\label{eq:tkm_grad_u_inc_quadrature}
    \nabla u_{\mathrm{inc}}(x)
    \approx
    \frac{\Delta k_x\Delta k_y\Delta k_z}{(2\pi)^3\epsilon_m}
    \sum_{\mu = 1}^{M_k}
    ik_\mu\,\widehat G_L(k_\mu)\widehat\rho_h(k_\mu)
    e^{ik_\mu\cdot x}.
\end{equation}
We use a 3D type-2
 non-uniform fast Fourier transform (NUFFT)~\cite{dutt1993fast,greengard2004accelerating,barnett2019parallel}
 to evaluate Eqs.~\eqref{eq:tkm_u_inc_quadrature} and \eqref{eq:tkm_grad_u_inc_quadrature} at a large, arbitrary set of target points simultaneously, at a cost that is quasi-linear in $M_k$ and number of target points. This is typically only one order of magnitude more expensive than an FFT.
 The inputs to the transform for the potential are the
 kernel-weighted Fourier coefficients $\widehat{\V u} = \{\widehat u(k_\mu) \}_{\mu = 1}^{M_k}$, where
\begin{equation}
\label{eq:tkm_coeff}
    \widehat u(k_\mu)
    =
    \widehat G_L(k_\mu)\widehat\rho_h(k_\mu).
\end{equation}
For the gradient components, the inputs are multiplied by the components of $ik_\mu$. The
stack of four type-2 NUFFTs can be done more efficiently in a single call.

\begin{remark}
  \label{r:skew}
  In practice, orbitals are often provided on regular but skew grids associated with the crystal lattice,
  so using a 3D FFT for \eqref{eq:rho_hat_discrete} would need resampling to a compatible Cartesian grid.
  Since there are typically $10^6$ such points, this is cumbersome.
  Instead we work directly with $\{y_j\}$ the input set of skew grid points,
  and $\omega_\rho$ their quadrature weight,
  evaluating \eqref{eq:rho_hat_discrete} using a single 3D type-1 NUFFT call. The
  speed difference vs an FFT is outweighed by the flexibility that comes from using NUFFTs for both input
  and target point sets.
\end{remark}

\subsection{Fourier-grid selection and error control}\label{sec:fourier_grid}

If $\Delta k_\alpha$ are chosen to obey certain bounds, then the aliasing error introduced by the discreteness of the Fourier quadrature in \eqref{eq:tkm_u_inc_quadrature} and \eqref{eq:tkm_grad_u_inc_quadrature} vanishes entirely, leaving only the high-frequency truncation error (due to the finite $k_{\max}$).

\begin{proposition}[Exactness of the infinite Fourier trapezoidal rule for the truncated kernel]
\label{prop:padding}
Let $\mathcal B$ have side lengths $(l_x,l_y,l_z)$, suppose
$\rho\in C_c^\infty(\mathbb R^3)$ with $\operatorname{supp}\rho\subset\mathcal B$, let
$h_{\mathrm{pad}}=c_{\mathrm{pad}}h$, and let $L$ be given by
Eq.~\eqref{eq:L_padded}. If
\begin{equation}
\label{eq:padding}
    \Delta k_\alpha
    \leq
    \frac{2\pi}{l_\alpha+h_{\mathrm{pad}}+L},
    \qquad \alpha\in\{x,y,z\},
\end{equation}
then the infinite-mode version of Eqs.~\eqref{eq:tkm_u_inc_quadrature} and \eqref{eq:tkm_grad_u_inc_quadrature} satisfies
\begin{equation}
\label{eq:tkm_infinite}
    \frac{\Delta k_x\Delta k_y\Delta k_z}
    {(2\pi)^3\epsilon_m}
    \sum_{\mu = 1}^{\infty}
    \widehat G_L(k_\mu)\widehat\rho(k_\mu)e^{ik_\mu\cdot x}
    =u_{\mathrm{inc}}(x),
    \qquad x\in\mathcal B_{\mathrm{pad}}
\end{equation}
and 
\begin{equation}
\label{eq:tkm_infinite_gradient}
    \frac{\Delta k_x\Delta k_y\Delta k_z}
    {(2\pi)^3\epsilon_m}
    \sum_{\mu = 1}^{\infty}
    ik_\mu\widehat G_L(k_\mu)\widehat\rho(k_\mu)e^{ik_\mu\cdot x}
    =\nabla u_{\mathrm{inc}}(x),
    \qquad x\in\mathcal B_{\mathrm{pad}},
\end{equation}
respectively, where $\sum_{\mu = 1}^{\infty}$ denotes $M_x,M_y,M_z\to\infty$ in Eq.~\eqref{eq:fourier_lattice}.
\end{proposition}

\begin{figure}[htbp]
    \centering
    \includegraphics[width=0.6\textwidth]{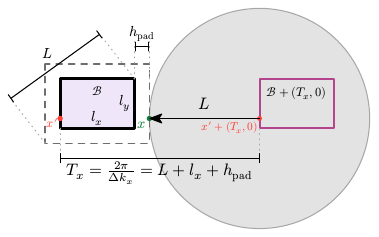}
    \caption{
        Two-dimensional illustration of the condition~\eqref{eq:padding}. 
        The source box $\mathcal B$ (black), padded target box $\mathcal B_{\mathrm{pad}}$ (dashed), and nearest periodic source image (purple) are shown at period $T_x=L+l_x+h_{\mathrm{pad}}$,
        in which case the radius-$L$ kernel support (gray disk) centered at one of the closest image source point just touches the target box.
    }
    \label{fig:tkm_partition}
\end{figure}

The proof relies on the finite spatial support of $G_L$, adapting arguments of \cite{VICO2016191} to the padded cuboidal case.
The 3D version of the Poisson summation formula~\cite[\S 7, Cor. 2.6]{SteinWeiss2016} gives
\begin{equation}
\label{eq:poisson_sum}
    \frac{\Delta k_x\Delta k_y\Delta k_z}{(2\pi)^3 \epsilon_m}
    \sum_{\mu = 1}^{\infty}
    \widehat G_L(k_\mu) \widehat \rho (k_\mu)e^{ik_\mu\cdot x}
    =
    u_{L}(x) + \sum_{n\in\mathbb Z^3\setminus\{0\}} u_{L}(x + T_n),
\end{equation}
where $T_n = (2\pi n_x/\Delta k_x, 2\pi n_y/\Delta k_y, 2\pi n_z/\Delta k_z)$, 
and the second term on the right-hand side is the error.
Since $\rho$ is supported in $\mathcal B$ and $G_L$ is supported in the ball of radius $L$, 
the support of $u_{L}$ lives within $\prod_{\alpha\in\{x,y,z\}} [- l_\alpha / 2 - L, l_\alpha / 2 + L]$.
Under condition~\eqref{eq:padding}, for all targets $x \in \mathcal B_{\mathrm{pad}}$ and $n \neq (0, 0, 0)$, each $x + T_n$ lies outside $u_{L}$'s support.
Thus the second term on the right-hand side of Eq.~\eqref{eq:poisson_sum} vanishes, and so does the error.
Recall that $u_{\mathrm{inc}}(x) = u_L(x)$ for all $x\in\mathcal B_{\mathrm{pad}}$, so that Eq.~\eqref{eq:tkm_infinite} follows.
Similar considerations apply to the gradient of $u_{L}$, leading to Eq.~\eqref{eq:tkm_infinite_gradient}.
Figure~\ref{fig:tkm_partition} illustrates the geometric separation that excludes contributions from periodic source images.

In practice, for efficiency we take equality in Eq.~\eqref{eq:padding}, giving the
largest $\Delta k_\alpha$ permitted by Proposition~\ref{prop:padding}, hence
the coarsest Fourier lattice satisfying the sufficient condition for eliminating aliasing error on $\mathcal B_{\mathrm{pad}}$.

\section{Efficient, parallel evaluation of the ERI tensor}\label{sec:assembly}

In this section, we describe the overall procedure for evaluating the ERI tensor in Eq.~\eqref{eq:four_indices} for a set of orbitals, which combines the BIE solver in Section~\ref{sec:bie} and the TKM/FMM solvers in Section~\ref{sec:tkm}.

\subsection{Batched BIE solves}\label{sec:blocksolve}

Evaluating the ERI tensor requires computing the potential generated by the orbital pair densities defined in Eq.~\eqref{eq:pair_density}.
For localized orbitals, many pair densities are small and can be neglected.
For a given threshold $\tau_{\mathrm{pair}}$, we define the set of \textit{active orbital pairs} as
\begin{equation}
\label{eq:retained_pairs}
    \mathcal I_{\mathrm{pair}}
    =
    \left\{
        (i,j):1\leq i\leq j\leq N_o,
        \ \norm{\rho_{ij}}_2>\tau_{\mathrm{pair}}
    \right\},
\end{equation}
where $\norm{\rho_{ij}}_2$ is the $L^2$ norm of the pair density $\rho_{ij}$.
Each $(i, j) \in \mathcal I_{\mathrm{pair}}$  defines an \textit{active density} $\rho_{ij}$ and a corresponding BIE to solve.
We enumerate the active densities as $\{\rho_a\}_{a=1}^{N_p}$, where $N_p=\abs{\mathcal I_{\mathrm{pair}}}$.

All active densities lead to BIEs with the same continuous operator but different right-hand sides and interface discretizations.
Instead of solving each of them separately, we group spatially nearby densities into batches so that their right-hand sides can be efficiently resolved on a common mesh.
The discrete operator~$\M A$ and FMM tree are shared across all source densities in a batch.
In practice, $\mathcal I_{\mathrm{pair}}$ is divided into $N_b$ batches $\{B_q\}_{q=1}^{N_b}$ with size $K_q\leq K_{\max}$, where $K_{\max}$ is a prescribed maximum batch size.

The geometry-driven criterion of Section~\ref{sec:edge_refinement} is source-independent, but the right-hand-side criterion of Section~\ref{sec:rhs_refinement} is not, since each right-hand-side $f_a$ requires refinement near its own source.
We therefore build the common mesh for a batch $\{\rho_a\}_{a = 1}^{K}$ using the following joint refinement criterion: 
\begin{equation}\label{eq:refine_criterion_batch}
    \max_{1\le a\le K} \norm{f_a-\Pi_P f_a}_{L^2(P)} <\tau_{\mathrm{rhs}}.
\end{equation}
Here we recall from \eqref{eq:rhs_interp} that $\Pi_{P} f_a$ is the interpolant of $f_a$.
This plays the role of \eqref{eq:refine_criterion}, i.e., we subdivide each panel until the criterion holds for all panels $P$.
Then the $K$ linear systems can be stacked as columns to give
\begin{equation}\label{eq:block_linsys}
    \M A \M \Sigma = \M F,
    \qquad
    \M \Sigma = [\V\sigma_1,\ldots,\V\sigma_K],
    \qquad
    \M F = [\V f_1,\ldots,\V f_K],
\end{equation}
with $\V f_a = -\partial_{\V n}u_{\mathrm{inc},a}$ sampled at the interface nodes.
The preconditioner \eqref{eq:linsys_precond} can then be applied.

Eq.~\eqref{eq:block_linsys} is solved using {\em block GMRES}~\cite{baker2006improving,parks2016block,onisk2023numerical}, which constructs a common block Krylov subspace for all $K$ systems, and is stopped at the selected residual tolerance $\tau_{\mathrm{GMRES}}$, as summarized in Algorithm~\ref{alg:solve}.
The FMM evaluates all $K$ strength vectors on the same source and target nodes, sharing the tree, interaction lists, and M2L matrix-vector multiplies, etc, across the batch.
The vectorized FMM has given several-fold acceleration in BEM applications~\cite{Makarov2021BEMFMM}.
The user-selected $K_{\max}$ controls the tradeoff between the FMM cost per right-hand side and the orthogonalization work and memory required by a wider block Krylov basis. 
We assess this tradeoff numerically in Section~\ref{sec:multirhs}.

\begin{algorithm}[H]
\footnotesize
\algrenewcommand\algorithmicindent{1.2em}
\algrenewcommand\algorithmicrequire{\textbf{Input:}}
\algrenewcommand\algorithmicensure{\textbf{Output:}}
\caption{Interface Poisson BIE solve for one source batch.}
\label{alg:solve}
    \begin{algorithmic}[1]
        \Require pair-density batch $B=\{\rho_a\}_{a=1}^{K}$ on uniform volume grids; interface $\Gamma$; permittivities $\{\epsilon_k\}$; order $p$; thresholds $\tau_{\mathrm{rhs}}$, $\tau_{\mathrm{quad}}$, $\tau_{\mathrm{GMRES}}$ and $\tau_{\rho}$; Fourier/NUFFT/FMM tolerance; smallest panel size $h_{\min}$
        \Ensure layer densities $\M\Sigma$ and the coefficient arrays $\{\widehat{\V u}_a\}$ of Eq.~\eqref{eq:tkm_coeff}
        \Procedure{SolveBatch}{$B$, $\Gamma$, $\{\epsilon_k\}$}
            \ForAll{$a = 1,\dots,K$}
                \State form $\widehat{\V u}_a$ of Eq.~\eqref{eq:tkm_coeff} from the input samples of $\rho_a$, using a type-1 NUFFT \Comment{Sec.~\ref{sec:tkm}}
            \EndFor
            \State partition $\Gamma$ into flat panels carrying $p{\times}p$ Gauss--Legendre rules; refine until the batch criterion Eq.~\eqref{eq:refine_criterion_batch} holds, then recursively subdivide near edges and corners until those panels have size $h_{\min}$ \Comment{Sec.~\ref{sec:adaptive_refinement}}
            \State construct $\M A = \M G + \M D^T + \M C^T$ of Eq.~\eqref{eq:opsplit} \Comment{Secs.~\ref{sec:quadrature},~\ref{sec:nearfield}}
            \State form RHS $\M F$ and solve $\M A\M\Sigma = \M F$ by block GMRES, applying $\M A$ via FMM plus the block-sparse correction
            \State \Return $\M\Sigma$, $\{\widehat{\V u}_a\}$
        \EndProcedure
    \end{algorithmic}
\end{algorithm}

\subsection{Post-evaluation of the potential}
\label{sec:eval}

For each active source density $\rho_a$, once the corresponding BIE has been solved, we evaluate the physical potential $\phi_a=u_{\mathrm{inc},a}+u_a$ on the volume grids~\eqref{eq:source_grid} of all active target densities $\{\rho_b\}_{b=1}^{N_p}$.
For each target density $\rho_b$, we then compute the corresponding ERI by the trapezoidal rule:
\begin{equation}\label{eq:V_contract}
    V_{ab} = \int_{\mathbb R^3} \rho_b(x) \phi_a(x)\,\mathrm{d}x \approx h^3 \sum_{m=1}^{N_{\rho_b}} \rho_b(y_m^{(b)}) \phi_a(y_m^{(b)}),
\end{equation}
where $\{y_m^{(b)}\}_{m=1}^{N_{\rho_b}}$ are the volume grid nodes of $\rho_b$.

More specifically, to compute $\{V_{ab}\}_{b = 1}^{N_p}$ for a given $a$,
we consider the union of {\em all} target grid points:
\begin{equation}
\label{eq:target_union}
    \mathcal{T}
    =
    \bigcup_{b=1}^{N_p}
    \left\{y_m^{(b)}:1\leq m\leq N_{\rho_b}\right\},
    \qquad
    N_{\mathcal T}=\abs{\mathcal T}.
\end{equation}
The incident potential $u_{\mathrm{inc},a}$ defined in Eq.~\eqref{eq:incident_potential_def} is then evaluated on $\mathcal T$ by the near/far procedure of Section~\ref{sec:tkm}.
The scattered potential $u_a$ is evaluated from the discretization of the single-layer representation \eqref{eq:single_layer}, using the density column vector $\V\sigma_a$ from the block linear solve \eqref{eq:block_linsys}.

Since the target set $\mathcal T$ can be close to the interface, special quadrature for near interactions is also needed.
For a panel $P$ and a target point $y$, the same criterion as Eq.~\eqref{eq:bernstein_criterion} is used to determine whether $y$ is near to or far from $P$.
Far-field contributions are evaluated by FMM-accelerated quadrature using the existing interface nodes and weights.
Integrals over near panels are evaluated by adaptive quadrature using the known density $\sigma_a$ to the desired accuracy, as a correction to the far-field contributions.
Unlike the on-surface correction, these near-field contributions are computed on-the-fly without forming or storing a correction matrix.
The resulting evaluation procedure is summarized in Algorithm~\ref{alg:eval}.

\begin{algorithm}[H]
\footnotesize
\algrenewcommand\algorithmicindent{1.2em}
\algrenewcommand\algorithmicrequire{\textbf{Input:}}
\algrenewcommand\algorithmicensure{\textbf{Output:}}
\caption{Evaluation of one source batch on the merged target set.}
\label{alg:eval}
    \begin{algorithmic}[1]
        \Require batch $B$ with global pair-index set $I_B$ and its solve panels and volume-lattice samples; layer densities $\M\Sigma$ and coefficient arrays $\{\widehat{\V u}_a\}$ from Algorithm~\ref{alg:solve}; merged set of all targets $\mathcal T$.
        \Ensure $\{V_{ab}\}_{b=1}^{N_p}$, $a \in I_B$
        \Procedure{EvalBatch}{$B$, $\M\Sigma$, $\{\widehat{\V u}_a\}$, $\mathcal T$}
            \ForAll{$a\in I_B$}
                \State evaluate $u_{\mathrm{inc},a}$ on $\mathcal T$ from $\widehat{\V u}_a$ and the samples of $\rho_a$ by the near--far split of Section~\ref{sec:tkm}
                \State evaluate $u_a$ on $\mathcal T$ from the discretized Eq.~\eqref{eq:single_layer}, with near panels selected by Eq.~\eqref{eq:bernstein_criterion}
                \State $\phi_a\gets u_{\mathrm{inc},a}+u_a$ on $\mathcal T$
                \For{$b=1,\dots,N_p$}
                    \State gather $\phi_a$ at $\{y_m^{(b)}\}$ and form $V_{ab}$ by Eq.~\eqref{eq:V_contract}
                \EndFor
            \EndFor
            \State \Return $\{V_{ab}\}_{b=1}^{N_p}$, $a \in I_B$
        \EndProcedure
    \end{algorithmic}
\end{algorithm}

\subsection{Distributed pipeline and computational cost}\label{sec:pipeline}

The overall procedure for evaluating the ERI tensor is summarized as a three-stage pipeline:
\begin{itemize}
    \item[1.] \emph{Prepare.} Sample and batch the pair densities, and construct target grids.
    \item[2.] \emph{Solve.} Construct the common discretization for each batch and solve Eq.~\eqref{eq:block_linsys}.
    \item[3.] \emph{Evaluate.} Evaluate each source potential to form the entries $V_{ab}$ with $b\geq a$.
\end{itemize}
The independent solve and evaluation stages are distributed across workers;
the full procedure is summarized in Algorithm~\ref{alg:campaign}.

\begin{algorithm}[H]
\footnotesize
\algrenewcommand\algorithmicrequire{\textbf{Input:}}
\algrenewcommand\algorithmicensure{\textbf{Output:}}
\caption{Distributed computation of the four-index tensor $\V V$.}
\label{alg:campaign}
    \begin{algorithmic}[1]
        \Require orbitals~$\V{\varphi}$; dielectric regions $\V \Omega$; pair-retention tolerance $\tau_{\mathrm{pair}}$; numerical controls $\tau_{\mathrm{rhs}}$, $\tau_{\mathrm{quad}}$, $\tau_{\mathrm{GMRES}}$, $\tau_{\rho}$, and the Fourier/NUFFT/FMM tolerance; maximum batch size~$K_{\max}$; smallest panel size $h_{\min}$; quadrature order $p$
        \Ensure non-negligible elements of the four-index tensor $\V V$ in sparse form
        \State set $\Gamma$; sample the pair densities on volume grids Eq.~\eqref{eq:source_grid}; select the pair densities using Eq.~\eqref{eq:retained_pairs}; form the $N_b$ localized batches of size at most $K_{\max}$; set $\mathcal{T}$ by Eq.~\eqref{eq:target_union} \Comment{serial}
        \ForAll{batches $B$ \textbf{in parallel}} \Comment{distributed across workers}
            \State $\M\Sigma_B, \{\widehat{\V u}_a\} \gets \textproc{SolveBatch}(B,\Gamma,\{\epsilon_k\})$ \Comment{Alg.~\ref{alg:solve}}
        \EndFor
        \ForAll{batches $B$ \textbf{in parallel}} \Comment{distributed across workers}
            \State $\{V_{ab}\} \gets \textproc{EvalBatch}(B,\M\Sigma_B,\{\widehat{\V u}_a\},\mathcal{T})$ \Comment{Alg.~\ref{alg:eval}}
        \EndFor
        \State \Return $\V V$
    \end{algorithmic}
\end{algorithm}

We estimate the main computational costs of Algorithm~\ref{alg:campaign}.
Let $N_p$ denote the number of active pair densities, $N$ a typical number of surface
Nystr\"om nodes, $N_\rho$ a typical number of volume-grid points per pair density, and $N_{\mathcal T}=\abs{\mathcal T}$ the number of distinct target points.
For localized orbitals, $N_p=\mathcal O(N_o)$.
The estimates below describe total computational work at fixed~$p$ and tolerances, assuming that the Fourier mode count $M_k = \mathcal O(N_\rho)$ and that local quadrature corrections only require linear work.
The solve and evaluation stages account for most of the measured runtime.

For the solve stage, evaluating the right-hand sides at the interface nodes costs
\begin{equation}
    C_{\mathrm f} = \mathcal O\bigl(N_p(N_\rho+N)\bigr) + \mathcal O\bigl(N_pN_\rho\log N_\rho\bigr),
\end{equation}
where the first term accounts for the FMM-accelerated far-field in Eq.~\eqref{eq:particle_gradient}, and the second term accounts for the NUFFT-based near-field evaluation in Eq.~\eqref{eq:tkm_grad_u_inc_quadrature}.
This is quasilinear in the total number of grid points $N_{\mathcal T}=\mathcal O(N_pN_\rho)$. 
The cost of BIE solves is dominated by applying $\M A$ in GMRES iterations, given by
\begin{equation}
    C_{\mathrm{BIE}} = \mathcal O\bigl(N_pN_{\mathrm{iter}}N\bigr),
\end{equation}
where $N_{\mathrm{iter}}$ is a representative number of block GMRES iterations, which is studied in Section~\ref{sec:convergence}.
For the evaluation stage, the scattered potentials are evaluated on $\mathcal T$, while the incident potentials are evaluated by the near/far procedure in Section~\ref{sec:tkm}, giving the cost
\begin{equation}
    C_{\mathrm{eval}} = \mathcal O\bigl( N_p(N + N_{\mathcal T}) \bigr) + \mathcal O\bigl(N_pN_\rho\log N_\rho\bigr).
\end{equation}
The first term is $\mathcal O(N_p^2N_\rho)$,
quadratic in the number of active pair densities (but still quasilinear with respect to number of grid points per pair).
This contrasts the previous stages, which were all $\bigO(N_p)$.

\section{Numerical Results}\label{sec:numerical_results}

In this section, we demonstrate the convergence, accuracy, timing and memory usage of the scheme.
We also show the method's performance on a large-scale, realistic system with complex dielectric geometry.
Our implementation is in Julia~\cite{doi:10.1137/141000671}, and is available at~\cite{biejl}.

\subsection{Convergence and effect of dielectric contrast}\label{sec:convergence}

We first study the convergence with respect to the quadrature order $p$ and the minimum panel size $h_\mathrm{min}$ for edge refinement at fixed dielectric contrast.
We then hold $p$ and $h_\mathrm{min}$ fixed and study the effect of increasing dielectric contrast on accuracy and GMRES iteration count, also showcasing the effect of the diagonal right preconditioner of Eq.~\eqref{eq:linsys_precond}.

\begin{figure}[htbp]
    \centering
    \includegraphics[width=0.99\textwidth]{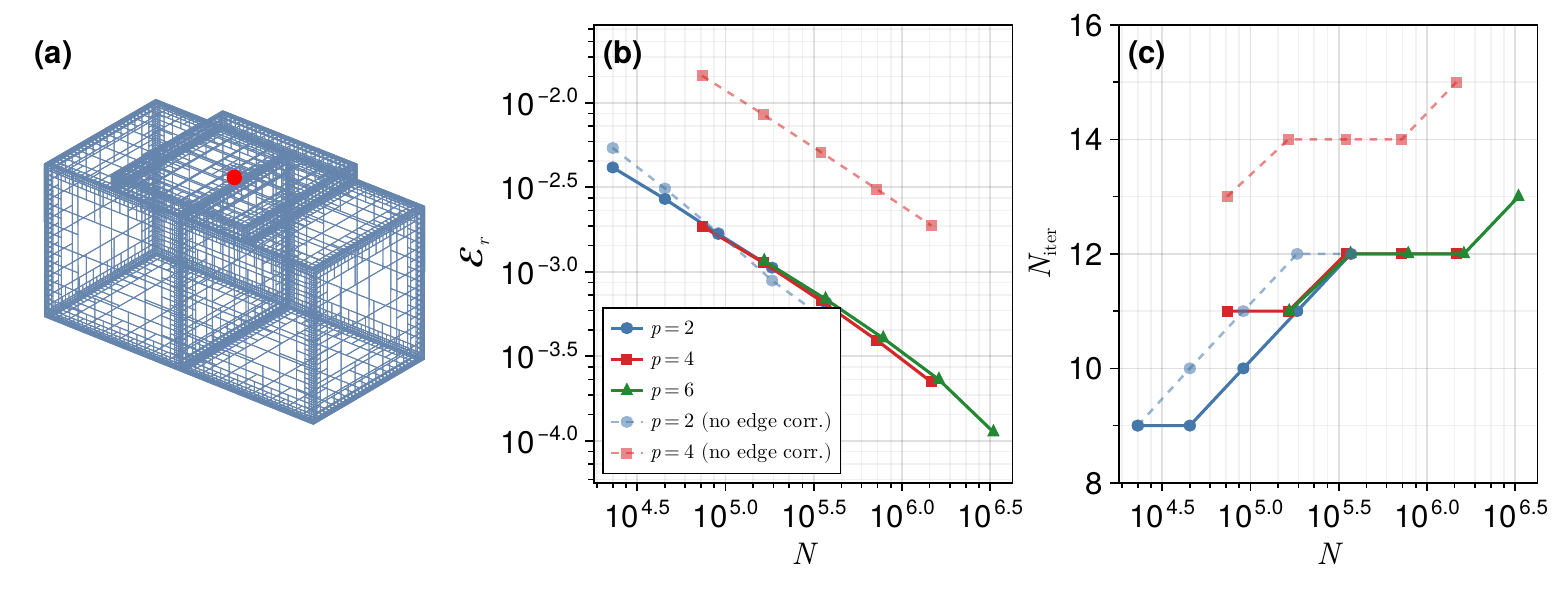}
    \caption{Three-dielectric benchmark with a single Gaussian source. (a) Geometry, showing adaptive interface panelization. (b) Relative error in the source self-interaction $U_0$ versus the number of Nystr\"om unknowns $N$.
    The label \textit{edge corr.} denotes the touching-panel correction described in Section~\ref{sec:nearfield}.
    (c) GMRES iteration count versus $N$.}
    \label{fig:fig1}
\end{figure}

A simple three-material geometry is used, as shown in Fig.~\ref{fig:fig1}(a): a $10 \times 10 \times 1$ dielectric slab ($\Omega_3$) rests on two $10\times10\times10$ cubes sharing a vertical face ($\Omega_1$ and $\Omega_2$), all embedded in vacuum.
A Gaussian density of width $0.05$ at the center of the slab is used as the source density~$\rho$, which is represented by the red dot in Fig.~\ref{fig:fig1}(a),
We compute $\rho$'s self-interaction $U_0$ as the quantity of interest for convergence studies, given by
\begin{equation}
    U_0 = \iint_{\mathcal{B} \times \mathcal{B}} \rho(x) W(x, x') \rho(x') \mathrm{d} x \mathrm{d} x',
\end{equation}
where $\mathcal{B}$ is the bounding box of $\rho$, and $W$ is the Green's function of the three-material system.

For the convergence study, we set $\epsilon_1=4$, $\epsilon_2=12$, and $\epsilon_3=10$, and decrease $h_{\mathrm{min}}$ from $1/2$ to $1/32$ for each of $p\in\{2,4,6\}$.
Other numerical tolerances are all set to $10^{-4}$, including the RHS-refinement, near-quadrature, GMRES, and volume-solver tolerances. 
The relative error 
\begin{equation}
    \mathcal E_r = \frac{\abs{U_0 - U_{\mathrm{ref}}}}{\abs{U_{\mathrm{ref}}}}
\end{equation}
is reported against the number of unknowns $N$. Here $U_{\mathrm{ref}}$ is from a reference calculation with $p = 8$, $h_{\mathrm{min}} = 1/64$ and other tolerances set to $10^{-6}$.
The corresponding errors and GMRES iteration counts $N_{\mathrm{iter}}$ are shown as solid curves in Fig.~\ref{fig:fig1}(b) and (c), respectively.
At comparable $N$, the error curves for different $p$ are almost identical, with convergence rate approximately $\mathcal{O}(N^{-0.7})$, suggesting that the error is dominated by resolving edge singularities.
Indeed, the analysis of Section~\ref{sec:edge_refinement} predicts $\mathcal{O}(h_{\min}^{\beta+1}) = \mathcal{O}(N^{-(\beta+1)})$, where $\beta+1\approx 0.72$ using the dominant $12{:}1$ dielectric ratio.
Meanwhile, $N_{\mathrm{iter}}$ remains between $9$ and $13$, growing only mildly under mesh refinement.
Dashed curves show results obtained without the touching-panel correction for $p=2$ and $p=4$.
At $p = 2$, the error and $N_{\mathrm{iter}}$ are similar to the corrected case, but at $p = 4$, the error is about an order of magnitude larger, and $N_{\mathrm{iter}}$ is larger.
Although $p=2$ yields comparable accuracy without the touching-panel correction, a small $p$ requires more extensive near-quadrature corrections, substantially increasing the computational cost.
This motivates using a larger $p$ together with touching-panel corrections to reduce errors associated with edge singularities.

We next isolate the effect of dielectric contrast, fixing $p=6$, $h_{\mathrm{min}}=1/32$, $\epsilon_1=4$, and $\epsilon_3=10$, and varying $\epsilon_2 \in \{6,20,60,200\}$.
As shown in Table~\ref{tab:contrast}, the relative error $\mathcal{E}_r$ increases mildly with $\epsilon_2$, from about $10^{-4}$ at $\epsilon_2=6$ to $10^{-3}$ at $\epsilon_2=200$.
The preconditioned GMRES iteration count also increases, from $11$ to $19$.
Thus increasing dielectric contrast increases the error and iteration count mildly at this fixed resolution.
We also see that the diagonal right preconditioner in Eq.~\eqref{eq:linsys_precond} reduces the iteration count by $39\%$--$56\%$, with the largest reduction at the highest contrast.

\begin{table}[th]
\centering
\caption{GMRES iteration counts $N_{\mathrm{iter}}$ and relative error $\mathcal{E}_r$ for different dielectric contrasts, with and without the diagonal right preconditioner of Eq.~\eqref{eq:linsys_precond}. Number of unknowns is about $3.3 \times 10^6$ for all cases.}
\label{tab:contrast}
\begin{tabular}{cccccc}
\hline
 & & & \multicolumn{2}{c}{$N_{\mathrm{iter}}$} & \\
\cline{4-5}
$\epsilon_1$ & $\epsilon_2$ & $\epsilon_3$ & unpreconditioned & preconditioned & $\mathcal{E}_r$ \\
\hline
$4$ & $6$   & $10$ & $18$ & $11$ & $1.1\times10^{-4}$ \\
$4$ & $20$  & $10$ & $28$ & $13$ & $1.7\times10^{-4}$ \\
$4$ & $60$  & $10$ & $36$ & $17$ & $3.5\times10^{-4}$ \\
$4$ & $200$ & $10$ & $43$ & $19$ & $1.0\times10^{-3}$ \\
\hline
\end{tabular}
\end{table}

\subsection{Performance of the batched interface solve}\label{sec:multirhs}

We benchmark the batched interface solve described in Section~\ref{sec:blocksolve},
focusing on strong scaling across CPU cores and the reduction in cost per source density obtained by solving $K$  right-hand sides simultaneously.
We report single-node timings for the surface discretization, BIE solve, and evaluate steps of Algorithms~\ref{alg:solve} and \ref{alg:eval}.
The discretization step includes the mesh discretization and the construction of the discrete operator $\M A$, and the BIE solve refers to the block GMRES solve of Eq.~\eqref{eq:block_linsys}.
Timings for the full distributed pipeline are presented in Section~\ref{sec:lattice}.

\begin{figure}[th]
    \centering
    \includegraphics[width=0.6\textwidth]{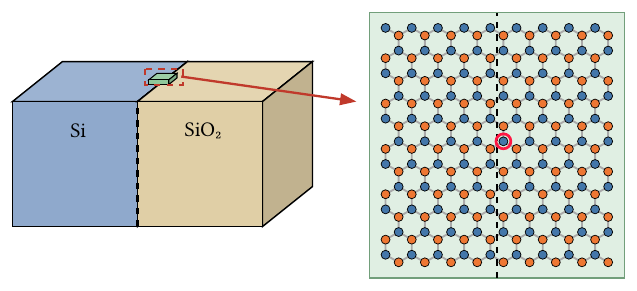}
    \caption{The triple-layer graphene above heterojunction system (left); top-down view of the graphene slab, with the $198$ carbon atoms colored by sublattice (right). The red circle marks the selected site associated with $\varphi_1$.}
    \label{fig:lattice_system}
\end{figure}

We consider a three-material system modeling
a trilayer of graphene placed above a Si/SiO$_2$ heterojunction, as shown in Fig.~\ref{fig:lattice_system}.
The heterojunction consists of two $270\times270\times270$~\AA{} cubes with relative permittivities of $11.9$ for Si and $3.9$ for SiO$_2$.
The graphene has thickness corresponding to three graphene layers, and is represented by a $43\times44\times9$~\AA{} slab with effective relative permittivity $\epsilon_{\mathrm{slab}}=2.4$, the cRPA value for graphene~\cite{PhysRevLett.106.236805}.

We focus orbitals in the central graphene layer, which contains $N_o = 198$ Wannier $p_z$ orbitals (see \ref{app:wannier}) centered on carbon atoms, as shown in the right panel of Fig.~\ref{fig:lattice_system}.
Each orbital density is supplied on a Wannier90~\cite{Pizzi2020} real-space skew
grid covering a $5\times5\times1$ graphene supercell with dimensions $12.325\times12.325\times15$~\AA{}.
The grid contains $150\times150\times192=4.32\times10^{6}$ points, with spacings of about $0.08$~\AA{}.
From those orbitals, we form the $N_p = 2613$ pair densities that are non-negligible according to the criterion of Eq.~\eqref{eq:retained_pairs}, with $\tau_{\mathrm{pair}} = 10^{-4}$.
This threshold corresponds approximately to only keep the pairs whose centers lie within $5$~\AA{} of one another.

Unless otherwise stated, all numerical tolerances are set to $10^{-5}$.
We use quadrature order $p=6$ and refine the interface mesh, and the edge panels are refined to a minimum size of $h_{\min}=1.14$~\AA{}. 
All benchmarks are performed on a single compute node equipped with two 48-core AMD EPYC 9474F processors and $1.5$~TB of RAM.

We first examine the strong scaling of the method for a single source density ($K=1$). 
We select the orbital closest to the center of the slab (see Fig.~\ref{fig:lattice_system}), denoted by $\varphi_1$, and use $\rho_{11}=\varphi_1^2$ as the source density. 
We solve the corresponding BIE, whose discretization requires $N\approx 1.35\times10^6$ interface unknowns. 
We then evaluate the resulting potential on the merged target set $\mathcal T$, which contains $7.17\times10^6$ points, and compute the $2613$ tensor elements $V_{ij11}$, where $(i, j) \in \mathcal{I}_{\text{pair}}$ (Eq.~\eqref{eq:retained_pairs}).
Keeping all other parameters fixed, we vary the number of threads from $1$ to $96$.

The results are shown in Table~\ref{tab:strongscaling}.
The total wall-clock time decreases from $1020$~s with $1$ thread to $62$~s with $96$ threads, corresponding to a speedup of $16.4\times$ and a parallel efficiency of $17\%$. 
The decrease in parallel efficiency is consistent with memory-bandwidth saturation during the FMM-dominated iterative solve and potential evaluation.

\begin{table}[th]
\centering
\caption{Strong scaling of Algorithm~\ref{alg:campaign} at batch size $K=1$ (one source density). The speedup and parallel efficiency are calculated relative to that of a single thread.}
\label{tab:strongscaling}
\begin{tabular}{ccccccc}
\hline
$N_{\mathrm{threads}}$ & discretization (s) & BIE solve (s) & eval (s) & total (s) & speedup & parallel efficiency \\
\hline
$1$ & $133.6$ & $759$ & $127$ & $1020$ & $1.0\times$ & $100\%$ \\
$4$ & $40.2$ & $381$ & $47$ & $469$ & $2.2\times$ & $54\%$ \\
$16$ & $14.4$ & $146$ & $22$ & $183$ & $5.6\times$ & $35\%$ \\
$64$ & $8.7$ & $50$ & $15$ & $74$ & $13.9\times$ & $22\%$ \\
$96$ & $8.9$ & $39$ & $14$ & $62$ & $16.4\times$ & $17\%$ \\
\hline
\end{tabular}
\end{table}

We then measure the benefit of batching multiple right-hand sides in the BIE solve and potential evaluation.
We fix the number of threads at $96$ and vary the batch size $K$ from $1$ to $79$,
by including more neighboring pair densities.
The results are shown in Table~\ref{tab:multirhs}.
The discretization time only increases slightly with $K$, from $10$~s to $16$~s.
Both the BIE solve and evaluation phases benefit from the batched FMM calls.
The end-to-end time per right-hand side decreases from $76$~s at $K=1$ to approximately $25$~s at $K=37$, and remains nearly constant through $K=58$, corresponding to a $2.5 \times$ increase in throughput.
Larger batches provide no further improvement, with the speedup declining slightly to $2.8\times$ at $K=79$.
Meanwhile, peak memory usage grows approximately linearly with $K$ from $24$~GB to about $1$~TB.
Thus, intermediate batch sizes, in this case specifically $K$ from $37$ to $58$, achieve nearly the highest measured throughput, with substantially lower memory usage than larger batches.

\begin{table}[th]
\centering
\caption{Wall-clock time and memory usage for the batched interface solve and evaluation at $N_{\mathrm{threads}}=96$ with different batch sizes $K$. The speedup is calculated relative to the total time per RHS at $K=1$.}
\label{tab:multirhs}
\begin{tabular}{cccccccc}
\hline
$K$ & $N_{\mathrm{iter}}$ & discretization (s) & BIE solve (s) & eval (s) & total (s) & peak RAM (GB) & speedup \\
\hline
$1$ & $11$ & $8.9$ & $39$ & $14$ & $62$ & $24$ & -- \\
$4$ & $10$ & $11.0$ & $88$ & $40$ & $140$ & $62$ & $1.8\times$ \\
$19$ & $10$ & $10.8$ & $317$ & $180$ & $507$ & $263$ & $2.3\times$ \\
$37$ & $9$ & $12.0$ & $552$ & $371$ & $934$ & $512$ & $2.5\times$ \\
$58$ & $9$ & $13.5$ & $846$ & $593$ & $1453$ & $783$ & $2.5\times$ \\
$79$ & $9$ & $15.8$ & $1258$ & $829$ & $2104$ & $1064$ & $2.3\times$ \\
\hline
\end{tabular}
\end{table}

\subsection{Large-scale ERI tensor computation}\label{sec:lattice}

We
apply the complete ERI computation pipeline, Algorithm~\ref{alg:campaign}, to the model introduced in Section~\ref{sec:multirhs}.
We first show the wall-clock time for the distributed computation of the four-index tensor $\V V$, and then show some representative entries of $\V V$ and their convergence properties.
We use the same numerical parameters: $p=6$, $h_{\min}=1.14$~\AA{} and the RHS-refinement, near-quadrature, and GMRES tolerances all set to $10^{-5}$.
The $N_o=198$ orbitals in the active graphene layer yield $N_p=2613$ pair densities.
From Table~\ref{tab:multirhs}, we select the batch size to be $K=46$ for the distributed computation, which gives $N_b=57$ batches in total.
Those batches are distributed across ten $96$-core nodes of the Rusty cluster at the Flatiron Institute,
and processed in parallel.

Table~\ref{tab:lattice_pipeline} shows the per-stage wall-clock time for the distributed computation of $\V V$.
The calculation takes $2$~h~$10$~m, or approximately $22$ node-hours (about $2100$ core-hours) in total. 
The solve and evaluation stages account for $50\%$ and $33\%$ of the wall-clock time, respectively.
We see that even with $N_p=2613$ pairs, the evaluation stage is not yet dominant over the other $\bigO(N_p)$-scaling stages despite its $\bigO(N_p^2)$ scaling.

\begin{table}[htbp]
\centering
\caption{Per-stage wall-clock time for the distributed computation of $\V V$. }
\label{tab:lattice_pipeline}
\begin{tabular}{lr}
\hline
stage & wall-clock \\
\hline
prepare      & $18$~m \\
solve        & $65$~m \\
eval         & $47$~m \\
\hline
total        & $130$~m \\
\hline
\end{tabular}
\end{table}

\begin{figure}[htbp]
    \centering
    \includegraphics[width=0.99\textwidth]{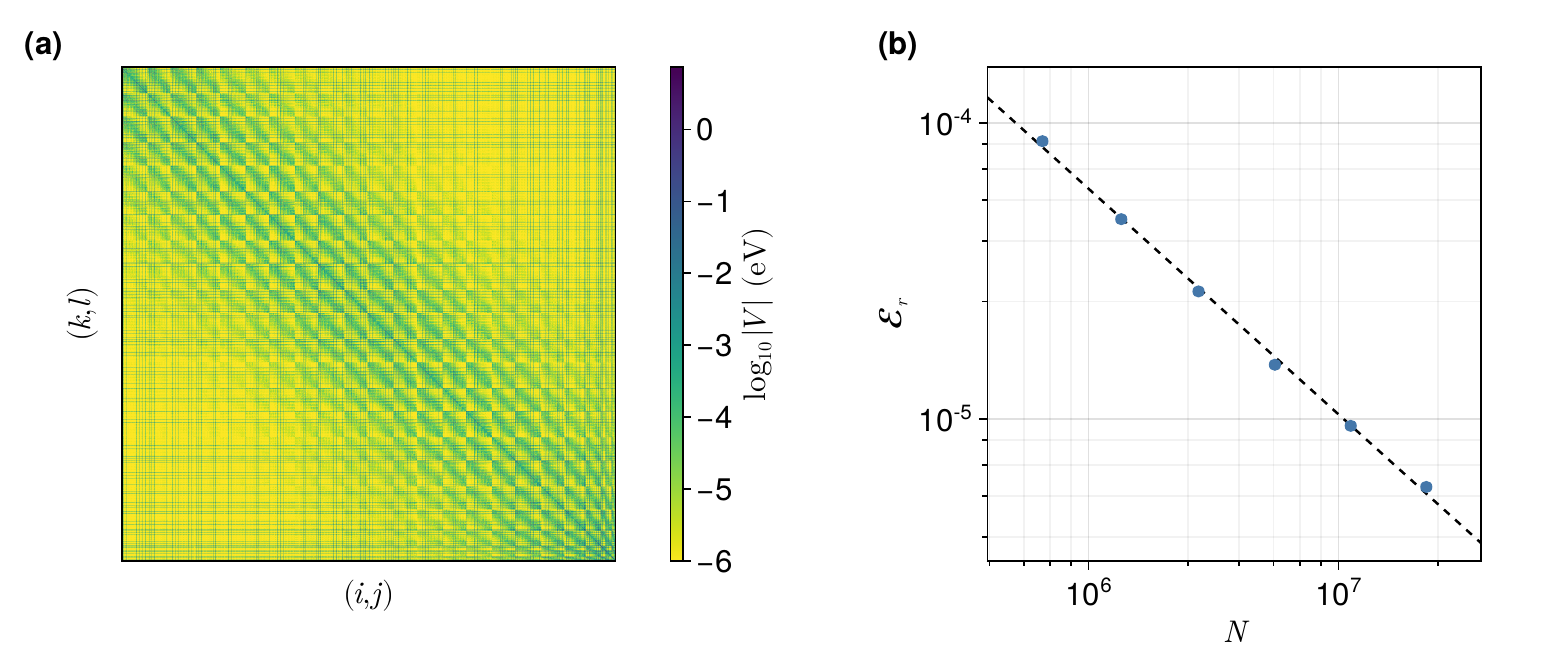}
    \caption{Results of the ERI tensor calculation.
    (a) $\V V$ in the locality-aware ordering, on a logarithmic scale in eV.
    (b) Self-convergence error of $U_1$ versus the number of interface unknowns $N$.
    The dashed line is a least-squares power-law fit.}
    \label{fig:lattice_V}
\end{figure}

Figure~\ref{fig:lattice_V}(a) shows the computed ERI tensor element magnitudes $|V_{ijkl}|$,
where $(i, j)$ and $(k, l)$ enumerate the $N_p$ orbital pairs.
Because the retained pairs are ordered according to spatial locality, the
matrix shows broad concentrations of relatively large entries together with
repeated horizontal and vertical structure.
To estimate the effect of discretization error at the production resolution, we perform a self-convergence study of $U_1:=V_{1111}$.
Holding all other parameters fixed, we decrease $h_{\min}$ to refine the interface mesh near edges and corners, increasing the number of interface unknowns from $6.5\times10^5$ to $2.2\times10^7$.
Figure~\ref{fig:lattice_V}(b) shows that the relative self-convergence error decays nearly as $N^{-0.8}$ over this range.
For $h_{\min}=1.14$~\AA{} and $N\approx1.35\times10^6$, the estimated relative error is approximately $4\times10^{-5}$ and the absolute error is about $0.3$~meV.

\begin{figure}[htbp]
    \centering
    \includegraphics[width=0.99\textwidth]{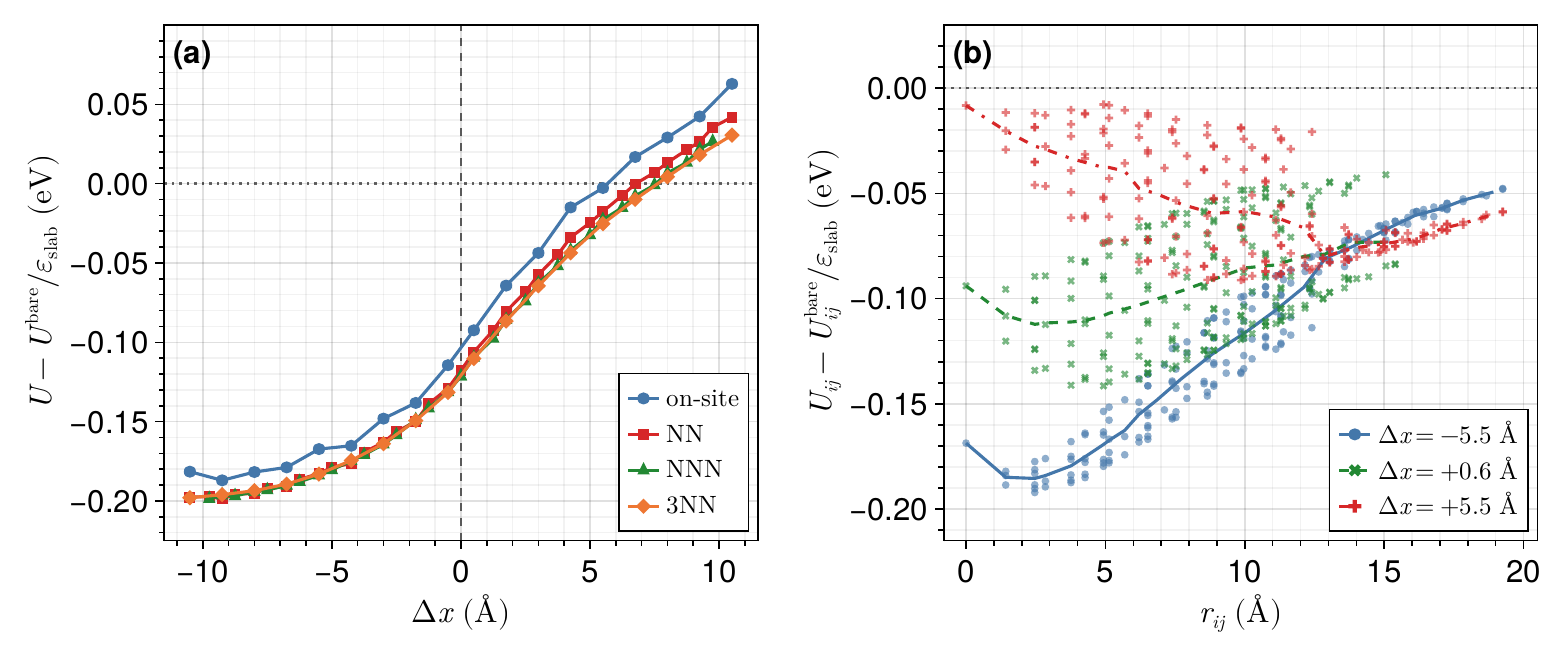}
    \caption{Effect of dielectric environment on the Coulomb interactions between orbitals.
        (a) On-site and first- to third-nearest-neighbor interaction shifts versus the signed distance $\Delta x$ from the junction (vertical dashed line), with Si on the left and SiO$_2$ on the right.
        (b) Inter-site interaction shifts versus orbital separation $r_{ij}$ for source sites at $\Delta x=-5.5$, $+0.6$, and $+5.5$~\AA{}.
        Points show individual matrix elements; curves show averages.
    }
    \label{fig:lattice_Ux}
\end{figure}

We next examine how the dielectric environment modifies the on-site interactions $U_i:=V_{iiii}$ and inter-site density--density interactions $U_{ij}:=V_{iijj}$.
Figure~\ref{fig:lattice_Ux} plots $\Delta U := U-U^{\mathrm{bare}}/\epsilon_{\mathrm{slab}}$, where $U^{\mathrm{bare}}$ is the bare Coulomb integral for the same orbitals.
In panel (a), the on-site and first- to third-nearest-neighbor interactions lie approximately $0.19$--$0.20$~eV below this reference on the Si side and increase by approximately $0.22$--$0.25$~eV across the junction, becoming slightly positive on the SiO$_2$ side.
Panel (b) compares $\Delta U_{ij}$ as function of orbital distance $r_{ij}$ for source orbitals $\varphi_i$ located in the Si region, near the junction, and in the SiO$_2$ region.
At short range, screening is strongest on the Si side, intermediate near the junction, and weakest on the SiO$_2$ side.
The results show that the Coulomb interactions are modified by the dielectric environment in a position-dependent manner, which explains the need for a full ERI tensor calculation in realistic device simulations.

\section{Conclusion}\label{sec:conclusion}

We have presented a boundary integral method and a scalable HPC implementation for computing the four-index electron repulsion integral tensor in piecewise-constant dielectric environments composed of cuboidal regions.
The method reformulates the governing Poisson interface problem as a second-kind boundary integral equation, discretized by a quad-panel-based Nystr\"om scheme, with $h$-adaptivity driven by both the incident data and the edge
singularities. %
To handle high aspect ratio cuboids efficiently, a Bernstein-radius criterion is used to engage
special quadrature corrections between relatively close or touching panels.
The incident fields of the smooth orbital-pair densities are evaluated by a truncated-kernel NUFFT method inside or near source supports, and quadrature farther away.
All far-field interactions are accelerated by a vectorized FMM, with block GMRES as the iterative solver,
showing significant amortization of costs over multiple right-hand sides.

Our extensive numerical studies included: i) convergence studies, showing rates matching an analytic density edge-singularity analysis,
ii) validation that GMRES iteration counts remain moderate under refinement or large dielectric contrasts,
iii) a strong scaling study with respect to thread count,
and iv) parallel timings for the full ERI tensor (6.9 million outputs) on $10^3$ cluster cores,
given 198 orbitals embedded in a realistic three-dielectric structure.
In this final set-up, the ratio of system length to orbital grid-spacing is about $10^4$.

Much of our focus was on the scalable many-right-hand-side workflow required for ERI construction.
While experiments confirmed that edge singularities induce a quite low convergence order
with respect to $N$ (lower than the expected $\bigO(N^{p/2})$ for uniform panels on a smooth surface),
an important conclusion is that 3--4 digit accuracy is achieved with $N\sim 10^6$, reasonable solution times,
and $h_{\min} \approx 1$\AA. This value of $h_{\min}$ has reached the atomic scale at which a sharp-edge model
becomes physically unimportant.
However, we plan to develop edge-singularity-capturing panel quadratures, which will allow us to increase the convergence order, and reduce $N$.
We also plan to extend the scheme to more general polyhedral or curved interfaces using existing BIE surface discretization tools.

For more general physics, possible future directions include periodic boundary conditions for extended multilayered materials, and the integration of the solver into self-consistent many-body frameworks, for which the screened interaction tensor computed here enters as a fundamental building block.

\appendix

\section{Density singularity exponent for a right-angle dielectric wedge}\label{app:wedge}

Here we derive the result on the density singularity
in Proposition~\ref{prop:edge_singularity},
by combining the jump relation with the right-angle special case of classical
dielectric-wedge analysis \cite{bobrovnikov73,andersen78}
(see also~\cite[\S5.3]{van2007electromagnetic}).
The generic singularity results from assuming invariance along the edge,
leaving a transverse 2D sectoral problem that is separable in polar coordinates $(r,\theta)$,
with the edge at the origin and the interfaces at $\theta = \pm \pi/4$.
The physical potential $\phi$
satisfies the 2D Laplace equation in each material,
\begin{equation*}
  \frac{\partial^2\phi}{\partial r^2}
  +\frac{1}{r}\frac{\partial\phi}{\partial r}
  +\frac{1}{r^2}\frac{\partial^2\phi}{\partial\theta^2}=0.
\end{equation*}
Substituting the separated-variable form $\phi(r,\theta)=R(r)\psi(\theta)$
with separation constant $g$ gives
\begin{equation*}
  r^2R''+rR'-g^2R=0,
  \qquad
  \psi''+g^2\psi=0.
\end{equation*}
For $g>0$, the bounded, nonconstant radial branch is $R(r)=r^g$.
Thus in each material sector $\phi(r,\theta)
= A r^g \cos(g\theta + C)$ for some coefficients $A,C$ that may differ between sectors.
To find eigenvalues $g>0$, we note that across either interface,
$\phi$ and $\epsilon \partial_{\V n} \phi$ are continuous.
Solutions fall into even and odd symmetry about $\theta=0$.
An even mode can be written
\begin{equation}
  \psi(\theta)=A\cos(g\theta), \quad \mbox{ for } |\theta|<\pi/4,
  \qquad \mbox{ and } \psi(\theta) = B\cos\bigl(g(\pi-\theta)\bigr) \quad \mbox{ otherwise}.
\end{equation}
Setting $a=\pi g/4$,
the two interface conditions at $\theta=\pi/4$ imply
\begin{equation}
  A\cos a=B\cos3a,
  \qquad -\epsilon_k A\sin a=\epsilon_{k'}B\sin3a.
\end{equation}
After using addition formulae, nontrivial solutions demand that
\begin{equation}\label{eq:right_angle_even_edge}
  \sin2a\,[2(\epsilon_k+\epsilon_{k'})\cos2a-(\epsilon_k-\epsilon_{k'})]=0.
\end{equation}
The sin factor vanishing gives only regular modes, while the other factor vanishing gives
\begin{equation}\label{eq:right_angle_even_exponent}
  g_e = \frac{4a}{\pi}
  =\frac{2}{\pi}
  \arccos\!\left(\frac{\epsilon_k-\epsilon_{k'}}{2(\epsilon_k+\epsilon_{k'})}\right)
    = \frac{2}{\pi} \arccos\!\left(\frac{1}{2 \gamma} \right),
\end{equation}
recalling \eqref{gammadef}, where the $\arccos$ is taken on the standard branch $[0,\pi]$.
Similarly, the smallest positive exponent for a odd mode 
is $g_o = \frac{2}{\pi} \arccos\!\left( - \frac{1}{2 \gamma} \right)$.
The smallest (leading) positive wedge exponent is therefore the smaller of the two,
\begin{equation}\label{eq:right_angle_leading_exponent}
  g_\ast = \min(g_e,g_o) = \frac{2}{\pi} \arccos\!\left(\frac{1}{2 \abs{\gamma}} \right),
\end{equation}
as in \cite{bobrovnikov73}.

Since $u_\text{inc}$ is smooth, \eqref{eq:potential_split} means that the scattered
potential $u$ has the same singularity as $\phi$.
Since the normal derivative is $\partial_{\V n} = r^{-1}\partial_\theta$,
then the radial dependence
of the jump in $\partial_{\V n} u$ is $\mathcal O(r^{g_\ast-1})$.
By the jump relation for the single-layer representation \eqref{eq:single_layer},
this jump determines $\sigma$, giving
Proposition~\ref{prop:edge_singularity} with $\beta = g_\ast-1$.
The condition $\beta>-1/3$ follows from $|\gamma|>1$.

If smooth potential
variation along the edge is allowed, one still gets the same singularity as the above
2D analysis, since a Fourier transform along the edge direction gives a 2D
modified Helmholtz equation in each material.
Then $R(r) = I_g(r)$, the regular modified Bessel function of real order $g$,
whose asymptotic power law is the same.
\section{Wannier function generation}\label{app:wannier}

The Wannier functions were obtained from density functional theory (DFT) calculations performed with Quantum ESPRESSO (v7.4) \cite{Giannozzi2009,Giannozzi2017} using the PBE exchange--correlation functional \cite{Perdew1996} and Grimme's DFT-D3 dispersion correction \cite{Grimme2010}.
Carbon was described by an optimized norm-conserving Vanderbilt pseudopotential from the SG15 library (C\_ONCV\_PBE-1.2) \cite{Hamann2013,Schlipf2015}, with a plane-wave kinetic-energy cutoff of 100~Ry.
We modeled monolayer graphene in a hexagonal unit cell with lattice constant $a = 2.465$~\AA{}, and periodic images were separated by $15$~\AA{} of vacuum.
The self-consistent ground state was converged to $10^{-10}$~Ry on a $\Gamma$-centered $32\times32\times1$ Monkhorst--Pack $k$-mesh, with $0.05$~eV ($3.7$~mRy) Gaussian smearing.
From these we built two maximally localized Wannier functions \cite{marzari1997maximally,marzari2012maximally} with C $p_z$ initial projections using Wannier90 (v3.1.0) \cite{Pizzi2020}.
The inner (frozen) energy window was set to $E_F \pm 3$~eV.
The Wannier functions were plotted on a $5\times5\times1$ real-space supercell for visualization.

\bibliographystyle{elsarticle-num-etal}
\bibliography{ref}

\end{document}